\documentclass[20pt,twoside,openany,reqno]{amsart}
 \usepackage{bbm}
 \usepackage{mathrsfs}
\input amssymb.sty
\usepackage{amsmath}
\usepackage[all,poly,knot]{xy}
\usepackage{amsthm}
\usepackage{amsfonts}
\usepackage{latexsym}
\usepackage{amscd}
\newcommand{\pf}{\noindent\begin {proof}}
\newcommand{\epf}{\end{proof}}

\def\bc{\begin{center}}
\def\ec{\end{center}}

\begin{document}

{\centerline{\LARGE {\bf On  cyclotomic elements and cyclotomic subgroups  }}

\bigskip

{\centerline{\LARGE {\bf in $K_{2}$ of a field}}

\bigskip

\bigskip

{\centerline{\it  Dedicated to Professor Jerzy Browkin}}

\bigskip

{\centerline{by}}

\bigskip

{\centerline{ { {\large K}EJIAN {\large X}U}\  \  AND \ {\large C}HAOCHAO {\large S}UN$^{\ast}$ \ }}

\begin{figure}[b]
\rule[-2.5truemm]{5cm}{0.1truemm}\\[2mm]
{\small
\begin{tabular}{ll}
$\ast$ This research is supported by National Natural Science
Foundation of China (No. 10871106).
\end{tabular}}
\end{figure}

\bigskip

\bigskip

{\centerline{\bf Abstract}}

\bigskip
The problem of expressing an element of $K_2(F)$ in a more explicit form gives rise to many works. To avoid a restrictive condition in a work of Tate, Browkin considered cyclotomic elements as the candidate for the element with an explicit form. In this paper, we modify and change Browkin's conjecture about cyclotomic elements into  more precise forms, in particular we introduce the conception of cyclotomic subgroup.
In the rational function field cases,  we determine completely the exact numbers of cyclotomic elements and cyclotomic subgroups contained in a subgroup generated by finitely many different cyclotomic elements; while in the number field cases, at first, a number field $F$ is constructed so that $G_5(F),$ the set of cyclotomic elements of order five in $K_2(F)$, contains at least three nontrivial cyclotomic subgroups, and then using Faltings' theorem on Mordell conjecture we prove that there exist  subgroups generated by an infinite number of cyclotomic elements to the power of some prime,  which contain no nontrivial cyclotomic elements.}

\bigskip


 \bigskip


\bigskip

\bigskip

{\centerline{\bf 1. Introduction}}

\bigskip

\noindent It follows from Matsumoto's theorem ([10]) that for a field $F,$ $K_2(F),$ the Milnor $K_2$-group,  can be generated by symbols $\{a, b\}, a, b \in F^*.$ In general, an element of $K_2(F)$ is only a product of symbols. Therefore, expressing an element of $K_2(F)$ in a simple and more explicit form  is highly expected. For a global field, Lenstra ([7]) proved a curious  fact that every element of $K_2(F)$ is not just a product of symbols, but  actually a symbol. More precisely, if $G$ is a finite subgroup of $K_2(F),$ then $G\subseteq \{a,F^*\}$ for some $a\in F^*.$

Furthermore,  for  a global field $F$ containing $\zeta_n,$ the $n$-th primitive root of unity, Tate ([18]) investigated the $n$-torsion of  $K_2(F)$ and  proved that
$$(K_2(F))_n=\{\zeta_n, F^{*}\}, \eqno(1. 1)$$
which implies that every element in the $n$-torsion $(K_2(F))_n$ can be written in the form of $\{\zeta_n, a\},$ where $a\in F^{*}.$
 Throughout this paper, for an abelian group $A,$ we  use the symbol $A_n$ to denote the $n$-torsion of $A,$ i.e., $A_n=\{a\in A |\, a^n=1\}.$
 Tate conjectured that the equality (1.1) is true for any field containing $\zeta_n.$ Mercurjev and Suslin proved Tate's conjecture ([9][18]).
Unfortunately, the condition $\zeta_n\in F$ is too restrictive. For example, as is well known, $K_2(\mathbb{Q})$ is a torsion group and contains elements of any order by Dirichlet's theorem.
But, by Tate's result, only elements of order $2$ in $K_2(\mathbb{Q})$ can be expressed explicitly.

In [1], Browkin considered  cyclotomic elements of $K_2(F),$ i.e., the elements of the form
$$c_n(a):=\{a, \Phi_{n}(a)\}, \ \ \ \   a, \Phi_n(a)\in F^{*},$$
where $\Phi_{n}(x)$ denotes the $n$-th cyclotomic polynomial.
The advantage of cyclotomic elements is that one can go without the condition $\zeta_n\in F.$

Let
 $$
G_{n}(F)=\{c_n(a) \in K_{2}(F)\mid {a, \Phi_{n}(a)} \in
F^{*}\}.$$
Then, Browkin ([1]) proved that $G_{n}(F)\subseteq(K_{2}(F))_{n} $,
i.e.,  $G_{n}(F)$ is contained in the $n$-torsion of
$K_2(F),$ in parcicular, he proved that for any field $F\neq \mathbb{F}_2,$ if $n=1,2,3,4,6$ and if $\zeta_n\in F,$ then every element  $\{\zeta_n, x\}\in K_2(F)$ can be written in the form $c_n(a).$
Moreover, it is also proved in [1]
 that $G_{n}(F)=(K_{2}(F))_{n} $ for $n=3$ and $F=\mathbb{Q}$ (for any field $F$ by Urbanowicz [20]).
As for $n=4,$ it follows from  [1]
for $F=\mathbb{Q}$ and Qin [12] for any field $F$ with ch$(F)\neq 2$ that every element of order 4 in $K_2(F)$ can be written in the form of
$c_4(a)\cdot v,$ where  $v\in K_2(F)$ with $v^2=1.$
But, in general,  as conjectured in [1], $G_n(F)$ is not a
group.

\bigskip

\noindent {\bf Browkin's Conjecture} ([1])\ {\it For any integer $n\neq 1,2,3,4$ or $6$
and any field $F$, $G_{n}(F)$ is not a subgroup of $K_{2}(F)$, in
particular, $G_{5}(\mathbb{Q})$ is not a subgroup of
$K_{2}(\mathbb{Q}).$ \vspace{1mm}}

\bigskip

Qin proved in [12] and [13] respectively that neither $G_5(\mathbb{Q})$ nor $G_7(\mathbb{Q})$ is a subgroup of $K_2(\mathbb{Q})$ and  that
$G_{2^{m}}(\mathbb{Q})$ is a group if and only if $n\leq 2.$ In [22], the authors proved that $G_{2^{n}3^{m}}(\mathbb{Q})$ is a group if and only if $n=2$ and $m=0$ (see [23] for more results).
But, the first author of the present paper  prove that for any number field $F,$ if $n\neq 4,8, 12$ is a positive integer having a square factor, then $G_{n}(F)$ is not a subgroup of $K_2(F)$ (see [23][25]). A similar result can be established for function fields ([25]).

 But, when $n$ is  a prime,  $G_n(F)$ seems difficult  to deal with in particular when $F$ is  a number field or, in general, a global field.
In [23], the authors have investigated the $l$-torsion of $K_2(F(x)),$ where $F(x)$ is the rational
function field over $F$ and $l$ is a prime with $l\neq $ch$(F),$ and proved that if $l\geq 5$  and if $\Phi_l(x)$ is irreducible in $F[x]$, then Browkin's  conjecture is true for
 $F(x)$. But we still do not know whether this is true for a number field.

Browkin's conjecture implies that corresponding to Tate's result, we could only expect results
on the ``outer structure'' of $G_n(F),$ that is, whether $(K_2(F))_n$ can be generated by something like $G_n(F).$
In fact, Lenstra proved that $(K_2(\mathbb{Q}))_5$
can be generated by $G_5(\mathbb{Q})$ (See [8][27]). In general, Qin conjectured that $(K_2(F))_n$ can be generated by all $G_m(F)$ with $ m$ being a divisor of $n,$ i.e.(See [2] for more results)

\bigskip

\noindent {\bf Qin's Conjecture} ([14]) \ $(K_2(F))_n=\langle G_m(F)\mid m|n\rangle.$

\bigskip

 Furthermore, the authors of [25] even conjectured that if $n$ has the factorization
$n=p_1^{e_1}p_2^{e_2}\cdots p_t^{e_t},$ then $(K_2(F))_n$ can be generated by all $G_{p_i^{m_i}}(F),$ i.e.,
$$(K_2(F))_n=\langle G_{p_i^{m_i}}(F)\mid 1\leq m_i\leq e_i, 1\leq i \leq t \rangle.$$

  In the present paper, we  turn to  the ``inner structure'' of $G_n(F)$, in particular, we are interested in the ``inner'' subgroup structure of $G_n(F)$. As a result,
  we   modify and change Browkin's conjecture into  more precise forms. A subgroup of $K_2(F)$ is called {\it cyclotomic} if  it is contained in $G_n(F)$.  Our questions are formulated  as follows.

\bigskip

\noindent {\bf Question 1} \ {\it   How many nontrivial cyclotomic elements  are there in
a subgroup of $K_2(F)$
 generated by finitely many different cyclotomic elements of order $n$ ?}

\noindent {\bf Question 2} \ {\it When $G_n(F)$ contains a nontrivial cyclotomic subgroup ?}

\noindent {\bf Question 3} \ {\it How many cyclotomic subgroups  are there in
a subgroup of $K_2(F)$
 generated by finitely many different cyclotomic elements of order $n$  ?}

\bigskip

It follows from [1] that for $F\neq \mathbb{F}_2$ and $n=1, 2,3, 4$ or $6,$ $G_n(F)$ itself is a cyclotomic subgroup  of $K_2(F).$
 In [28]and [29], the authors proved that for a local field $F$, $G_n(F)$ is a cyclotomic subgroup in most cases (see also [4]). Moreover, they conjectured that for a local field $F,$ $G_n(F)$ is always a cyclotomic subgroup of $K_2(F)$. But for  a number field, the picture seems  different. From [26], we  only know that  a subgroup of $K_2(F(x))$ generated by a cyclotomic element
contains at least two  non-cyclotomic elements.

 In this paper, for the rational function field $F(x),$ we give more precise result, that is,
 we determine the exact number of nontrivial cyclotomic elements and nontrivial cyclotomic subgroups in a subgroup
  generated by some cyclotomic  elements in $G_l(F(x))\subseteq K_2(F(x)),$  where $l$ is a prime with $l\neq$ch$(F).$ More precisely,
let $\mathfrak{G}_{l}(n; F)$ denote a subgroup of $K_2(F(x))$ generated by $n$ essentially
  distinct (see section 4) cyclotomic elements of some kind in $G_l(F(x)),$ and let $c(\mathfrak{G}_{l}(n; F))$ and $cs(\mathfrak{G}_{l}(n; F))$ denote respectively the numbers of nontrivial cyclotomic elements and nontrivial cyclotomic subgroups contained in  $\mathfrak{G}_{l}(n; F),$ then we prove the following result (See Theorem 5.17).

\bigskip

 {\bf Theorem 1.1}\ \, {\it Assume that $l\geq 5$ is a prime number and  $F$ is a field
such that $\Phi_l(x)$ is irreducible in $F[x]$. Let $n$ be a positive integer satisfying
$$ n\leq \frac{l-3}{2}.$$

 i) If ch$(F)=0, $ then $c(\mathfrak{G}_{l}(n; F))=2n,$  and so $cs(\mathfrak{G}_{l}(n; F))=0.$

 ii) If ch$(F)=p\neq0,$ then $c(\mathfrak{G}_{l}(n; F))=n(2+|\mathfrak{Z}(l,p)|),$  where
$$\mathfrak{Z}(l,p):=\{t \mid 2\leq t \leq l-2, \ t\equiv p^{2m} \mbox{or} -p^{2m} (\mbox{mod}\, l)\ \mbox{for some}\ m\in \mathbb{N} \}.$$

iii) If ch$(F)=p\neq0,$ then we have
  $$cs(\mathfrak{G}_{l}(n; F))> 0 \Longleftrightarrow  \ l\equiv 3\, (\mbox{mod}\, 4) \ \mbox{and}
\ p \mbox{ is a primitive root of} \ l.$$
   In this case,
  $cs(\mathfrak{G}_{l}(n; F))=n,$ i.e.,  $\mathfrak{G}_{l}(n; F)$ contains exactly $n$ nontrivial  cyclotomic subgroups.

iv) Every nontrivial cyclotomic subgroup of $\mathfrak{G}_{l}(n; F)$ is a cyclic subgroup of order $l,$ i.e.,
every nontrivial cyclotomic subgroup has the form $\mathfrak{G}_{l}(1; F).$}

\bigskip

   We do not know how to remove the condition $n\leq \frac{l-3}{2}$ in  Theorem 1.1.   We present some computations for the cases
   $n> \frac{l-3}{2},$ in particular for $n=2,3.$ The results of computations coincide
   with the above theorem. So it seems that the condition $n\leq \frac{l-3}{2}$ is removable.

As for the number field cases, the situation seems quite different. In the proof
of Theorem 1.1, the essential use is made of the fact that the
function field $F(x)$ has a nontrivial derivation. Thus the proof
does not carry over to number fields.

However, it seems curious that we can really construct a number field $F$ and a cyclotomic element in $K_2(F)$ such that the cube of this element is also cyclotomic (we can do some things for the square), and  as a consequence,  we can  construct a number field $F$ so that  $G_5(F)$ contains  a
nontrivial  cyclotomic subgroup. Furthermore, we can also construct a number field $F$ so that  $G_5(F)$ contains  at least three nontrivial  cyclotomic subgroups !

But, it seems that this is not true in general. In fact,  using Faltings' theorem on Mordell conjecture, we can prove the following result (see Theorem 10.4).

\bigskip

{\bf Theorem 1.2} \ \ {\it Assume that $F$ is a number field and
 $n\neq 1,4,8,12$ is a positive integer. If there is a prime $p$ such that  $p^2|n,$
then  there exist infinitely many nontrivial cyclotomic elements $\alpha_1, \alpha_2, \ldots, \alpha_m, \ldots  \in
G_{n}(F)$ so that
 $$\langle\alpha_1^{p}\rangle\subsetneq \langle\alpha_1^{p}, \alpha_2^p\rangle\subsetneq \ldots \subsetneq \langle\alpha_1^{p}, \alpha_2^p, \ldots, \alpha_m^p\rangle\subsetneq \ldots $$
 and
$$\langle\alpha_1^{p}, \alpha_2^p, \ldots, \alpha_m^p, \ldots \ \rangle\cap G_{n}(F)=\{1\}.$$
}

This implies that in $K_2(F)$ there exists a subgroup generated by cyclotomic elements to the power of some prime, which  contains no nontrivial cyclotomic elements.
Clearly, this result is  more precise than Browkin's conjecture.  Hence, in general,  for a number field $F,$ we conjecture that if $p> 5$ is a prime,
then $G_p(F)$ contains no nontrivial cyclotomic subgroups.

This paper is organized as follows. The first part of this paper, i.e. from  section 2 to section 8,  focuses  on the case of function fields. In section 2, we discuss some basic properties relative to cyclotomic polynomials; in section 3, the definition of tame homomorphism and its computation are given;
in section 4, to remove superfluous generators in a finitely generated subgroup of  $K_2(F(x)),$ we introduce the conception of `essentially distinct elements'; while in section 5, our aim is to prove  Theorem 1.1. In section 6,  some computations are presented for the case $n>\frac{l-3}{2},$ in particular, for $n=2$ or $3;$ in section 7, for the preparation  of the next section, two diophantine equations are discussed; while in section 8, a further example is given. Then, in the second part of this paper, we consider the number field cases. More precisely, in section 9,  we construct the cube (resp. square) of
some cyclotomic element which is also cyclotomic and as a result some cyclotomic subgroups of order five are constructed, in particular   a number field $F$ is constructed so that  $G_5(F)$ contains  at least three nontrivial  cyclotomic subgroups, and  finally in section 10,  Theorem 1.2 is proved.

\bigskip

\bigskip

{\centerline{\bf 2. Cyclotomic Polynomials.
}}

\bigskip

\noindent  Let $l\geq 5$ be a prime number and  $F$ a field of characteristic $\neq l.$ Through out this paper we will always assume that the cyclotomic polynomial $\Phi_l(x)$ is irreducible in $F[x].$ We denote by $\zeta$ any root of $\Phi_l(x).$

Let $\Phi_l(x,y):=y^{l-1}\Phi_l(x/y).$ From the irreducibility of $\Phi_l(x)$ in $F[x]$ it follows the irreducibility of $\Phi_l(x,y)$ in $F[x,y].$

\bigskip

\noindent  {\bf Theorem 2.1}\ \ {\it For any nonzero polynomial $f(x), g(x)\in F[x]$ we have
  }
  $$\mbox{deg}\Phi_l(f(x),g(x))=(l-1)\cdot \mbox{max}(\mbox{deg}f(x),\mbox{deg}g(x)).$$

\noindent  {\it Proof:} \ We have
 $$\Phi_l(f(x),g(x))=f(x)^{l-1}+f(x)^{l-2}g(x)+\ldots + g(x)^{l-1}.\eqno(2.1)$$
Let $a_0x^r$ and $b_0x^s$ be leading terms of polynomials of $f(x)$ and $g(x),$ respectively.

 If $r\neq s,$ say, $r> s,$ then, by (2.1), the leading term of $\Phi_l(f(x),g(x))$  is equal to $(a_0x^r)^{l-1}=a_0^{l-1}x^{(l-1)r}.$

If $r=s,$ then all summands in (2.1) are of same degree, and the sum of their leading terms is
$$(a_0x^r)^{l-1}+(a_0x^r)^{l-2}b_0x^r+\ldots +(b_0x^r)^{l-1}=\Phi_l(a_0,b_0)x^{r(l-1)}.$$
Moreover, $\Phi_l(a_0,b_0)=b_0^{l-1}\Phi_l(a_0/b_0)\neq 0,$ since the polynomial $\Phi_l(x)$  irreducible in $F[x]$ can not have a zero in $F.$

 Thus in both cases the leading term of $\Phi_l(f(x),g(x))$ is of degree $(l-1)r=(l-1)\cdot$max(deg$f(x),$ deg$g(x)$). \hfill $\Box$

\bigskip

\noindent {\bf Theorem 2.2}\ \  {\it If $f(x), g(x)\in F[x]$ are relatively prime polynomials, then the degree of every factor of $\Phi_l(f(x), g(x))$ is
divisible by $l-1.$}

\noindent  {\it Proof:} It is sufficient to prove that the degree of every irreducible factor of $\Phi_l(f(x), g(x))$ is divisible by $l-1.$

In $F[x]$ we have
$$\Phi_l(f(x), g(x))=\prod_{j=1}^{l-1}(f(x)-\zeta^{j}g(x)).\eqno(2.2)$$
Let $\alpha$ be a root of an irreducible factor $h(x)$ of $\Phi_l(f(x), g(x)).$ Then it is a root of $\Phi_l(f(x), g(x)),$ hence, by (2.2), $f(\alpha)-\zeta^jg(\alpha)=0$ for some $1\leq j\leq l-1.$

 Therefore $f(\alpha)=0$ if and only if $g(\alpha)=0.$ It follows that $f(\alpha)g(\alpha)\neq 0,$ since $f(x)$ and $g(x)$ cannot have a common root, because they are relatively prime.

 Consequently $\zeta^j=f(\alpha)/g(\alpha)\in F(\alpha).$ Hence $F(\zeta)\subseteq F(\alpha).$ Therefore
 $$\mbox{deg}h(x)=(F(\alpha):F)=(F(\alpha):F(\zeta))(F(\zeta):F)=(F(\alpha):F(\zeta))(l-1),$$
 since $\alpha$ and $\zeta$ are roots of polynomials $h(x)$ and $\Phi_l(x),$ respectively, which are irreducible.

 \hfill$\Box$

\bigskip

\noindent  {\bf Corollary 2.3} \ {\it If }max(deg{\it $f(x),$}deg{\it $g(x))=1,$ then the polynomial $\Phi_l(f(x), g(x))$ is irreducible.}

\noindent {\it Proof:} By Theorem 2.1, deg$\Phi_l(f(x),g(x))=l-1,$ and, by Theorem 2.2, every factor of $\Phi_l(f(x),g(x))$ has degree divisible by $l-1.$
Therefore the polynomial $\Phi_l(f(x),g(x))$ has only one factor, so it is irreducible. \hfill$\Box$

\bigskip

\noindent {\bf Theorem 2.4} \ {\it Let $f(x), g(x)\in F[x]$ satisfy $(f(x), g(x))=1$ and} deg{\it $f(x)\geq 1.$ Let $\mathfrak{p}$ be the ideal of $F[x]$ generated by an irreducible factor of $\Phi_l(f(x),g(x)).$

Then for $r\in \mathbb{Z}$
$$(f(x)/g(x))^r\equiv 1 (\mbox{mod}\ \mathfrak{p})\ \ \mbox{if and only if}\ \ l \mid r.$$ }

\noindent  {\it Proof:} \ Since $\mathfrak{p}$ is generated by an irreducible polynomial, it is a prime ideal of $F[x].$
 From $\Phi_l(f(x),g(x))\mid f^l-g^l,$ it follows that $f^l\equiv g^l ($mod $ \mathfrak{p}),$ and $g(x)\not\equiv 0($mod $\mathfrak{p}),$
 because $f(x)$ and $g(x)$ are relatively prime. Hence $(f(x)/g(x))^{l}\equiv 1($mod $\mathfrak{p}).$

 If $l\nmid r$ and $(f(x)/g(x))^{r}\equiv 1($mod $\mathfrak{p}),$ then from the last two congruences it follows that $f(x)/g(x)\equiv 1($mod $\mathfrak{p}),$ i.e. $f(x)\equiv g(x)$(mod $\mathfrak{p}).$ Hence
 $$\Phi_l(f(x),g(x))=\sum_{j=0}^{l-1}f(x)^jg(x)^{(l-1)-j}\equiv lg(x)^{l-1}(\mbox{mod}\, \mathfrak{p}),$$
 then $g(x)\equiv 0 (\mbox{mod}\, \mathfrak{p}),$ which is impossible.
The contradiction shows that $l|r.$

 Conversely, if $l|r,$ then from the congruence $(f(x)/g(x))^l\equiv 1 $(mod $\mathfrak{p}$), it follows that
 $(f(x)/g(x))^r\equiv 1 $(mod $\mathfrak{p}$).
\hfill  $\Box$

\bigskip

 Let $W(F)$ be the group of roots of unity in $F.$

 We say that matrices $A,B\in GL(2,F) $ are {\it essentially distinct} if
 $$B\neq \alpha \left (\begin{matrix}
   \mu & 0\\
   0 & 1
  \end{matrix}\right )\left (\begin{matrix}
   0 & 1\\
   1 & 0
  \end{matrix}\right )^{\epsilon}A$$
 for every $\alpha\in F^{*}, \mu \in W(F),$ and $ \epsilon=0$ or $1.$

 Thus if $A=\left (\begin{matrix}
   a & b\\
   c & d
  \end{matrix}\right )\in GL(2, F)$ then all matrices which are not essentially distinct from $A$ are
  $$\alpha \left (\begin{matrix}
   \mu a & \mu b\\
   c & d
  \end{matrix}\right )\ \  \mbox{and}\ \ \alpha \left (\begin{matrix}
   \mu c & \mu d\\
   a & b
  \end{matrix}\right ), \ \mbox{for all}\ \alpha\in F^*, \mu \in W(F).$$

\bigskip

\noindent {\bf Theorem 2.5} \  {\it If matrices
$$\left (\begin{matrix}
    a_1 & b_1\\
   c_1 & d_1
  \end{matrix}\right ),\  \ \left (\begin{matrix}
    a_2 & b_2\\
   c_2 & d_2
  \end{matrix}\right )\in GL(2,F)$$
  are essentially distinct, then the polynomials
  $$\Phi_l(a_1x+b_1,c_1x+d_1)\ \ \mbox{and}\ \ \Phi_l(a_2x+b_2,c_2x+d_2)$$
  are relatively prime.
}

\noindent {\it Proof:} If matrices $A_1$ and $A_2$ are essentially distinct, then for every $B\in GL(2,F)$ the matrices $A_1B$ and $A_2B$ are essentially distinct. Therefore, taking $B=A_1^{-1}$ we can assume that $A_1=I$ is the identity matrix, and $A_2= \left (\begin{matrix}
    a_2 & b_2\\
   c_2 & d_2
  \end{matrix}\right ).$

Assume that the corresponding polynomials $\Phi_l(x)$ and $\Phi_l(ax+b,cx+d)$ are not relatively prime. Since they are irreducible and of the
same degree, they differ by a constant factor only :
$$\Phi_l(x)=\alpha\Phi_l(ax+b,cx+d)\ \ \mbox{for some}\ \alpha\in F^*.$$
Hence corresponding linear factors of both polynomials differ by a constant factor, in particular
$$x-\zeta=\alpha_1((ax+b)-\zeta^r(cx+d))\ \ \mbox{for some}\ \alpha_1\in F(\zeta)^*\ \mbox{and}\ 1\leq r \leq l-1.$$
Comparing coefficients we get
$$1=\alpha_1(a-\zeta^rc),\ \ \ -\zeta=\alpha_1(b-\zeta^rd).$$
Eliminating $\alpha_1$ we obtain
$$-\zeta(a-\zeta^rc)=b-\zeta^rd.$$

If $r\neq 1,l-1,$ then $1,\zeta, \zeta^r,\zeta^{r+1}$ are linearly independent over $F,$ hence $a=b=c=d=0,$ which is impossible.

If $r=1,$ then $-\zeta a+\zeta^2 c=b-\zeta d$ implies that $b=c=0$ and $a=d.$ Consequently $A_2=a\left (\begin{matrix}
    1 & 0\\
   0 & 1
  \end{matrix}\right )$ is not essentially distinct from $A_1=I.$

If $r=l-1,$ then $-\zeta a+c=b-\zeta^{l-1}d$ implies that $a=d=0, b=c.$ Consequently $A_2=b\left (\begin{matrix}
    0 & 1\\
   1 & 0
  \end{matrix}\right )$ is also not essentially distinct from $A_1=I.$

  In every case we get a contradiction. Therefore the polynomials $\Phi_l(x)$ and $\Phi_l(ax+b,cx+d)$ are relatively prime. \hfill $\Box$

\bigskip

\bigskip

{\centerline{\bf 3. Tame Homomorphisms
}}

\bigskip

\noindent For a nonzero prime ideal $\mathfrak{p}$ of $F[x],$ the tame homomorphism
$$\tau_{\mathfrak{p}}:\ K_2(F(x))\longrightarrow (F[x]/\mathfrak{p})^{*}$$
is defined by
$$
\tau_{\mathfrak{p}}(\{u,v\})\equiv (-1)^{v_{\mathfrak{p}}(u)v_{\mathfrak{p}}(v)} \frac{u^{v_{\mathfrak{p}}(v)}}{v^{v_{\mathfrak{p}}(u)}}
(\text{mod}\, \mathfrak{p}), \eqno(3.1)$$
where $u,v\in F(x)^{*}.$

\bigskip

\noindent {\bf Lemma 3.1} \ {\it Let $f(x), g(x)\in F[x]$ satisfy $(f(x), g(x))=1,$ \mbox{deg}$f(x)g(x)> 0.$ For a nonzero prime ideal $\mathfrak{p}$ of $F[x]$ denote $r_{\mathfrak{p}}:=v_{\mathfrak{p}}(\Phi_l(f(x),g(x))).$}

(i) {\it We have
$$\tau_{\mathfrak{p}}\Big(c_l\Big(\frac{f}{g}\Big)\Big)\equiv \left\{\begin{array}{ll}
                (\frac{f}{g})^{r_{\mathfrak{p}}}\not\equiv 1 (\bmod \mathfrak{p}), & \  \text{if} \ l\nmid r_{\mathfrak{p}},\\
                1 (\bmod \mathfrak{p}), & \   \text{if} \ l\mid r_{\mathfrak{p}}.
            \end{array}
        \right.
$$}

(ii) {\it In particular, if}\, max(deg$f(x),$ deg$g(x))=1,$ {\it then
$$\tau_{\mathfrak{p}}\Big(c_l\Big(\frac{f}{g}\Big)\Big)\equiv \left\{\begin{array}{ll}
                \frac{f}{g} \not\equiv 1 (\bmod \mathfrak{p}), & \  \text{if} \ \mathfrak{p}=(\Phi_l(f(x),g(x))),\\
                1 (\bmod \mathfrak{p}), & \ \text{otherwise}.
            \end{array}
        \right.
$$
}

\noindent {\it Proof:} (i) From $(f(x),g(x))=1$ it follows that $(f(x)g(x), \Phi_l(f(x),g(x)))=1.$ Therefore for every prime ideal $\mathfrak{p}$ of $F[x]$
at most one of the numbers $v_{\mathfrak{p}}(f(x)), v_{\mathfrak{p}}(g(x)),$ $ v_{\mathfrak{p}}(\Phi_l(f(x),g(x))$ does not vanish.

Clearly, we have
$$c_l\Big(\frac{f(x)}{g(x)}\Big)=\Big\{\frac{f(x)}{g(x)}, \Phi_l\Big(\frac{f(x)}{g(x)}\Big)\Big\}=\Big\{\frac{f(x)}{g(x)}, \Phi_l(f(x),g(x))\Big\}\{f(x),g(x)\}^{-(l-1)},\eqno(3.2)$$
because $\{g(x), g(x)^2\}=1$ and $l-1$ is even.

If $v_{\mathfrak{p}}(f(x))> 0$ and $r_{\mathfrak{p}}=0,$ then $\Phi_l(f(x),g(x))\equiv g(x)^{l-1}(\mbox{mod} \, \mathfrak{p}).$ Hence, by (3.1) and (3.2),
$$\tau_{\mathfrak{p}}\Big(c_l\Big(\frac{f}{g}\Big)\Big)\equiv \Phi_l(f(x),g(x))^{-v_{\mathfrak{p}}(f(x))}g(x)^{(l-1)v_{\mathfrak{p}}(f(x))}\equiv 1 (\mbox{mod}\, \mathfrak{p}).$$

If $v_{\mathfrak{p}}(g(x))> 0$  and $r_{\mathfrak{p}}=0,$  then we prove similarly that $\tau_{\mathfrak{p}}(c_l(f(x)/g(x)))\equiv 1 (\mbox{mod}\, \mathfrak{p}).$

If $v_{\mathfrak{p}}(f(x))=v_{\mathfrak{p}}(g(x))=0$ and $r_{\mathfrak{p}}=0,$ then (3.2) implies that
$\tau_{\mathfrak{p}}(c_l(f(x)/g(x)))\equiv 1(\mbox{mod}\, \mathfrak{p}).$

If $r_{\mathfrak{p}}> 0,$ then, by (3.1) and (3.2), $\tau_{\mathfrak{p}}(c_l(f(x)/g(x)))\equiv (f(x)/g(x))^{r_{\mathfrak{p}}}(\mbox{mod}\, \mathfrak{p}).$

Moreover, by Theorem 2.4, $(f(x)/g(x))^{r_{\mathfrak{p}}}\not\equiv 1 (\mbox{mod}\mathfrak{p})$ if and only if $l\nmid r_{\mathfrak{p}}.$

(ii) By Corollary 2.3, the polynomial $\Phi_l(f(x),g(x))$
is irreducible. Therefore  $r_{\mathfrak{p}}=v_{\mathfrak{p}}(\Phi_l(f(x),g(x)))=1.$ It is sufficient to apply the first part of the theorem with $r_{\mathfrak{p}}=1.$

\hfill $\Box$

\bigskip

\bigskip

{\centerline{\bf 4. Essentially Distinct Elements
}}

\bigskip

\noindent  It is well-known that
$$PGL(2,F):=GL(2,F)/ Z,$$
where $Z$ is the center of $GL(2,F)$, that is,
 $Z=F^{*}\cdot \left (\begin{matrix}
   1 & 0\\
   0 & 1
  \end{matrix}\right ).$
  Similarly,
  $$PSL(2,F):=SL(2,F)/Z \subset PGL(2,F).$$

  In the following, we will use the symbol $\overline{\left (\begin{matrix}
   a & b\\
   c & d
  \end{matrix}\right )}$ to denote  the images of  $\left (\begin{matrix}
   a & b\\
   c & d
  \end{matrix}\right )$ in $PGL(2, F).$ Clearly, the element $c_l\big(\frac{ax+b}{cx+d}\big)$ depends only on the coset $\overline{\left (\begin{matrix}
   a & b\\
   c & d
  \end{matrix}\right )}.$

\bigskip

We will focus on   the following subset of $G_l(F(x)):$
$$GG_l(F(x)):=\Big \{c_l\Big(\frac{ax+b}{cx+d}\Big)\in G_l(F(x))\mid \overline{\left (\begin{matrix}
   a & b\\
   c & d
  \end{matrix}\right )}\in PGL(2,F)\Big\}.$$
 $$SG_l(F(x)):=\left\{c_l\Big(\frac{ax+b}{cx+d}\Big)\in GG_l(F(x))\mid \left (\begin{matrix}
   a & b\\
   c & d
  \end{matrix}\right )\in SL(2,F)\right\},$$
  $$TG_l(F(x)):=\left\{c_l(x+b)\in GG_l(F(x))\ | \  b\in F \right\}.$$

\bigskip

\noindent   {\bf Definition 4.1} \ {\it Let
 $$\alpha=c_l\Big(\frac{a_1x+b_1}{c_1x+d_1}\Big), \ \  \beta = c_l\Big(\frac{a_2x+b_2}{c_2x+d_2}\Big)\in GG_l(F(x)).$$
We say that  $\alpha, \beta$  are {\it essentially distinct}
   if the matrices $\left (\begin{matrix}
   a_1 & b_1\\
   c_1 & d_1
  \end{matrix}\right )$ and $\left (\begin{matrix}
   a_2 & b_2\\
   c_2 & d_2
  \end{matrix}\right )$ are essentially distinct.}

\bigskip

\noindent  {\bf Lemma 4.2}\ \ {\it Assume that $\Phi_l(x)$ is irreducible in $F[x].$ Let
 $$\alpha=c_l\Big(\frac{a_1x+b_1}{c_1x+d_1}\Big), \ \  \beta = c_l\Big(\frac{a_2x+b_2}{c_2x+d_2}\Big)\in  GG_l(F(x)).$$

  If $\alpha=\beta,$ then we have
$$\frac{a_1x+b_1}{c_1x+d_1}=\frac{a_2x+b_2}{c_2x+d_2}.$$
}

\noindent  {\it Proof:}
Since $\Phi_{l}(x)$ is irreducible in $F[x],$ so are $\Phi_l(a_ix+b_i, \, c_ix+d_i)) \, (i=1,2)$ from Corollary 2.3.  From Lemma 3.1 we have the following congruences:
$$\tau_{\mathfrak{p}}(\alpha)\equiv \left\{\begin{array}{ll}
                \frac{a_1x+b_1}{c_1x+d_1}\not\equiv 1\ (\bmod \mathfrak{p}), & \  \text{if} \ \mathfrak{p}=(\Phi_{l}(a_1x+b_1, \, c_1x+d_1)),\\
                1 (\bmod \mathfrak{p}), & \ \text{otherwise};
            \end{array}
        \right.$$
$$\tau_{\mathfrak{p}}(\beta) \equiv \left\{\begin{array}{ll}
                \frac{a_2x+b_2}{c_2x+d_2}\not\equiv 1\ (\bmod \mathfrak{p}), & \ \text{if} \ \mathfrak{p}=(\Phi_{l}(a_2x+b_2, \, c_2x+d_2)),\\
                1 (\bmod \mathfrak{p}), & \ \text{otherwise}.
            \end{array}
        \right.$$

\bigskip

 If $\alpha =\beta,$ then we have $\tau_{\mathfrak{p}}(\alpha)=\tau_{\mathfrak{p}}(\beta),$  so we must have $(\Phi_{l}(a_1x+b_1, \, c_1x+d_1))= (\Phi_{l}(a_2x+b_2, \, c_2x+d_2))$ as primes. Therefore,  for the prime $\mathfrak{p}=(\Phi_{l}(a_1x+b_1, \, c_1x+d_1))$ we have
$$\frac{a_1x+b_1}{c_1x+d_1}\equiv \frac{a_2x+b_2}{c_2x+d_2} (\bmod \mathfrak{p}).$$
So
$$(a_1x+b_1)(c_2x+d_2)=(a_2x+b_2)(c_1x+d_1),$$
that is,
$$\frac{a_1x+b_1}{c_1x+d_1}=\frac{a_2x+b_2}{c_2x+d_2},$$
as required. \hfill $\Box$

\bigskip

\noindent  {\bf Lemma 4.3}\ \ {\it Assume that $\Phi_l(x)$ is irreducible in $F[x].$ Let
 $$\alpha=c_l\Big(\frac{a_1x+b_1}{c_1x+d_1}\Big), \ \  \beta = c_l\Big(\frac{a_2x+b_2}{c_2x+d_2}\Big)\in  GG_l(F(x)).$$
Then
$$\alpha=\beta \Longleftrightarrow \overline{\left (\begin{matrix}
   a_1 & b_1\\
   c_1 & d_1
  \end{matrix}\right )}=\overline{\left (\begin{matrix}
   a_2 & b_2\\
   c_2 & d_2
  \end{matrix}\right )}\in PGL(2,F).$$

In particular, if  $\alpha, \beta \in GG_l(F(x))$   are essentially distinct,  then they must be different, i.e., $\alpha\neq \beta.$ }

\noindent  {\it Proof:} $\Leftarrow:$ Clear.

$\Rightarrow:$ \
If $\alpha =\beta,$ then from Lemma 4.2 we have
 $$\frac{a_1x+b_1}{c_1x+d_1}=\frac{a_2x+b_2}{c_2x+d_2}.$$
So
$$a_1c_2=a_2c_1,\ \ b_1d_2=b_2d_1,$$
$$a_1d_2+b_1c_2=a_2d_1+b_2c_1.$$

\bigskip

Assume that $a_1c_2=a_2c_1=0.$ If $a_1=0,$ then $b_1c_1\neq 0$ since $a_1d_1-b_1c_1\neq 0,$ so $a_2=0,$ therefore $b_2c_2\neq0 $ since $a_2d_2-b_2c_2\neq 0.$
So, we can let $\frac{d_1}{d_2}=\frac{b_1}{b_2}=\frac{c_1}{c_2}=u\neq 0.$
Then  we have
$$\left (\begin{matrix}
   0 & b_1\\
   c_1 & d_1
  \end{matrix}\right )=u\left (\begin{matrix}
   0 & b_2\\
   c_2 & d_2
  \end{matrix}\right ),$$
  so
  $$\overline{\left (\begin{matrix}
   0 & b_1\\
   c_1 & d_1
  \end{matrix}\right )}=\overline{\left (\begin{matrix}
   0 & b_2\\
   c_2 & d_2
  \end{matrix}\right )}.$$

If  $c_2=0, $ the result is the same. Therefore, we should have $a_1c_2=a_2c_1\neq0.$ Similarly,  we have $b_1d_2=b_2d_1\neq0.$

\bigskip

Let $\frac{a_1}{a_2}=\frac{c_1}{c_2}=u\neq 0$ and $\frac{b_1}{b_2}=\frac{d_1}{d_2}=v\neq 0.$ Then from $a_1d_2+b_1c_2=a_2d_1+b_2c_1,$
we have $(a_2d_2-b_2c_2)(u-v)=0,$ which leads to $u=v$ since $a_2d_2-b_2c_2\neq 0.$ Hence, we get
$$\left (\begin{matrix}
   a_1 & b_1\\
   c_1 & d_1
  \end{matrix}\right )=u\left (\begin{matrix}
   a_2 & b_2\\
   c_2 & d_2
  \end{matrix}\right ),$$
so
$$\overline{\left (\begin{matrix}
   a_1 & b_1\\
   c_1 & d_1
  \end{matrix}\right )}=\overline{\left (\begin{matrix}
   a_2 & b_2\\
   c_2 & d_2
  \end{matrix}\right )}.$$

The last part of the lemma is obvious.  \hfill$\Box$

\bigskip

\noindent  {\bf Corollary 4.4}\ \ {\it Let
 $$\alpha=c_l(x+b_1),  \ \  \beta =c_l(x+b_2).$$
 Then the following statements are equivalent.}

 i) \ {\it $\alpha$ and $\beta$ are essentially distinct.}

 ii) \ {\it $\alpha\neq\beta.$}

 iii) \  {\it $b_1\neq  b_2.$}
 }

\noindent {\it Proof:} \ i)$\Rightarrow$ ii) \ It follows from Lemma 4.3.

ii)$\Rightarrow$ iii) \ Clear.

iii)$\Rightarrow $i) \ It is easy to check directly that $\left (\begin{matrix}
   1 & b_1\\
   0 & 1
  \end{matrix}\right )$ essentially distinct from  $\left (\begin{matrix}
   1 & b_2 \\
   0& 1
  \end{matrix}\right )$
 if and only if
 $b_1\neq  b_2.$ \hfill$\Box$

\bigskip

In the following, we will use the symbols $\mathfrak{G}_{l}(n; F), \mathfrak{S}_l(n;F)$ and $\mathfrak{T}_l(n;F)$ to denote  subgroups of $K_2(F(x))$
generated by (any)  $n$ essentially distinct nontrivial elements in $GG_l(F(x)), $ $SG_l(F(x))$ and $TG_l(F(x)),$ respectively.

From Corollary 4.4, we have

\bigskip

\noindent {\bf Lemma 4.5} \ {\it There exist some mutually different $b_1, b_2, \ldots,b_n \in F$  such that $\mathfrak{T}_l(n; F)$ is generated as follows.
$$\mathfrak{T}_l(n; F)=\langle c_l(x+b_1),  \ldots, c_l(x+b_n) \rangle.$$
}
\hfill$\Box$

 A subgroup of $K_2(F)$ is called {\it cyclotomic} if  it is contained in $G_n(F)$. In general, for a subgroup $H$ of $K_2(F(x)),$ we use the symbol $c(H)$ (resp. $cs(H)$)to denote the number of cyclotomic element(resp. cyclotomic subgroup )of $H$

\bigskip

\bigskip

{\centerline{\bf 5. The Rational Function Field Case
}}

\bigskip

 \noindent Assume that $l\geq 5$ is a prime. Let
$$\beta =\prod_{i=1}^{n} c_l\Big(\frac{a_ix+b_i}{c_ix+d_i}\Big)^{l_i},\eqno(5.1)$$
 where
 $1\leq l_i \leq l-1$ and $n\geq 1.$  If $n\geq 2,$ we assume that
$$\left (\begin{matrix}
   a_i & b_i\\
   c_i & d_i
  \end{matrix}\right )\in GL(2,F), \ \ 1\leq i \leq n,$$
  are essentially different from each other.

  It is well known that
  $$Gal(F(x)/F)\cong PGL(2, F)$$
and that $PGL(2,F)$ acts as automorphisms on $K_2(F(x))$ through
  $$\sigma\cdot \{f(x), g(x)\}:=\{f(x), g(x)\}^{\sigma}=\{f(\sigma (x)), g(\sigma (x))\}=\Big\{f\Big(\frac{ax+b}{cx+d}\Big), g\Big(\frac{ax+b}{cx+d}\Big)\Big\},$$
  where $\sigma=\overline{\left (\begin{matrix}
   a & b\\
   c & d
  \end{matrix}\right )} \in PGL(2,F)$ with $\sigma(x)=\frac{ax+b}{cx+d}.$

Applying the automorphism of the field $F(x),$ we may assume that the first factor on the right hand side of (5.1) is $c_l(x)^{l_1}.$

The polynomials $\Phi_l(a_ix+b_i,c_ix+d_i)$ are irreducible and by Theorem 2.5, pairwise relatively prime, hence the ideals $\mathfrak{p}_i:=(\Phi_l(a_ix+b_i,c_ix+d_i))$ in $F[x]$ for $i=1,2,\ldots, n$ are prime and  distinct.

We will prove some necessary conditions for $\beta$ to be cyclotomic. First we investigate the factorization of $\Phi_l(f(x),g(x)).$

\bigskip

\noindent {\bf Theorem 5.1}\ {\it Assume that the element $\beta$ given by} (5.1) {\it is cyclotomic:
$$\beta =c_l \Big (\frac{f(x)}{g(x)}\Big), \eqno(5.2)$$
where $f(x), g(x)\in F[x], (f(x),g(x))=1,$} deg{\it$(f(x)g(x))\geq 1.$ }

(i) {\it Then
$$\Phi_l(f(x), g(x))=\alpha \Psi^l \prod_{i=1}^{n}\Phi_l(a_ix+b_i, c_ix+d_i)^{r_i},\eqno(5.3)$$
where $\alpha \in F^*, \Psi\in F[x],$ and $ r_i:=v_{\mathfrak{p}_i}(\Phi_l(f(x), g(x)))$ satisfies $l \nmid r_i.$ We have $l-1|$}deg$\Psi.$

(ii) {\it Moreover,
$$\Big(\frac{f(x)}{g(x)}\Big)^{r_i}\equiv \Big(\frac{a_ix+b_i}{c_ix+d_i}\Big)^{l_i}\not\equiv 1\ (\mbox{mod}\, \mathfrak{p}_i)\ \mbox{for}\ i=1,2,\ldots,n. \eqno(5.4)$$
}

\noindent {\it Proof:} By Lemma 3.1 (i), for every prime ideal $\mathfrak{p}$ of $F[x]$ we have
$$\tau_{\mathfrak{p}}\Big(c_l\Big(\frac{f(x)}{g(x)}\Big)\Big)\equiv \left\{\begin{array}{ll}
                (\frac{f(x)}{g(x)})^{r_{\mathfrak{p}}}\not\equiv 1 (\bmod\, \mathfrak{p}), & \  \text{if} \ l\nmid r_{\mathfrak{p}},\\
                1 (\bmod\, \mathfrak{p}), & \   \text{if} \ l\mid r_{\mathfrak{p}},
            \end{array}
        \right.\eqno(5.5)
$$
and, by Lemma 3.1 (ii),
$$\tau_{\mathfrak{p}}\Big(c_l\Big(\frac{a_ix+b_i}{c_ix+d_i}\Big)\Big)\equiv \left\{\begin{array}{ll}
                \frac{a_ix+b_i}{c_ix+d_i}\not\equiv 1 (\bmod\, \mathfrak{p}), & \  \text{if} \ \mathfrak{p}=\mathfrak{p}_i,\\
                1 (\bmod\, \mathfrak{p}), & \   \mbox{otherwise.}
            \end{array}
        \right.\eqno(5.6)
$$
From (5.1) and (5.2) we get
$$c_l \Big (\frac{f(x)}{g(x)}\Big)=\prod_{i=1}^{n} c_l\Big(\frac{a_ix+b_i}{c_ix+d_i}\Big)^{l_i},\eqno(5.7)$$
where $1\leq l_i \leq l-1.$

Applying the tame homomorphism $\tau_{\mathfrak{p}},$ where $\mathfrak{p}$ is any prime ideal of $F[x],$ to both sides of (5.7), in view of (5.5) and (5.6), we obtain
$$\tau_{\mathfrak{p}}\Big(c_l\Big(\frac{f(x)}{g(x)}\Big)\Big)\not\equiv 1 (\mbox{mod}\, \mathfrak{p}) \Longleftrightarrow l\nmid v_{\mathfrak{p}}(\Phi_l(f(x),g(x)))\Longleftrightarrow \mathfrak{p}\in \{\mathfrak{p}_1,\mathfrak{p}_2,\ldots,\mathfrak{p}_n\}. \eqno(5.8)$$

It follows that in the representation of $\Phi_l(f(x),g(x))$ as the product of powers of relatively prime polynomials, the irreducible factors $\Phi_l(a_ix+b_i,c_ix+d_i)$ appear with the exponents $r_i$ not divisible by $l,$ and other factors appear with the exponents divisible by $l.$ This proves (5.3).

The divisibility $l-1|$deg$\Psi$ follows from Theorem 2.2, since $\Psi$ is a factor of $\Phi_l(f(x),g(x)).$ Thus we have proved (i).

By (5.8), $\tau_{\mathfrak{p}_i}\Big(c_l\Big(\frac{f(x)}{g(x)}\Big)\Big)\equiv \Big(\frac{f(x)}{g(x)}\Big)^{r_i}\not\equiv 1  (\mbox{mod}\, \mathfrak{p}_i)$ and, by (5.6),
$$\tau_{\mathfrak{p}_i}\Big(c_l\Big(\frac{a_jx+b_j}{c_jx+d_j}\Big)\Big)\equiv \left\{\begin{array}{ll}
                \frac{a_ix+b_i}{c_ix+d_i}\not\equiv 1 (\bmod\, \mathfrak{p}_i), & \  \text{if} \ j=i,\\
                1 (\bmod\, \mathfrak{p}_i), & \   \text{if} \ j\neq i.
            \end{array}
        \right.
$$
Consequently (5.7) implies that
$$\Big(\frac{f(x)}{g(x)}\Big)^{r_i}\equiv \Big(\frac{a_ix+b_i}{c_ix+d_i}\Big)^{l_i}\not\equiv 1\ (\mbox{mod}\, \mathfrak{p}_i),$$
which proves (ii). \hfill$\Box$

\bigskip

\noindent {\bf Theorem 5.2}\ {\it Denote $\theta:=$}max(deg$f(x),$deg$g(x)$). {\it Moreover, let $n\geq 2.$ Under the assumption of Theorem 5.1 we have
$$l\leq 2\theta+1.$$}
\noindent {\it Proof:} By  Theorem 2.4, we have
$$\Big(\frac{f(x)}{g(x)}\Big)^l\equiv \Big(\frac{a_ix+b_i}{c_ix+d_i}\Big)^{l}\equiv 1 (\mbox{mod}\, \mathfrak{p}_i).$$
Therefore raising both sides of (5.4) to the exponent $r^{\prime}_{i}$ such that $r_ir^{\prime}_i\equiv 1$(mod $l$), we get
$$\frac{f(x)}{g(x)}\equiv \Big(\frac{a_ix+b_i}{c_ix+d_i}\Big)^{m_i} (\mbox{mod} \, \mathfrak{p}_i),$$
where $1\leq m_i \leq l-1, m_i \equiv l_i r_i^{\prime}$ (mod $l$). Hence
$$\frac{f(x)}{g(x)}\equiv \Big(\frac{c_ix+d_i}{a_ix+b_i}\Big)^{l-m_i} (\mbox{mod}\, \mathfrak{p}_i),$$
From $\mathfrak{p}_i=(\Phi_l(a_ix+b_i, c_ix+d_i))$ we deduce that
$$
\begin{array}{ll}
\Phi_l(a_ix+b_i, c_ix+d_i)\mid f(x)(c_ix+d_i)^{m_i}-g(x)(a_ix+b_i)^{m_i},\\ \Phi_l(a_ix+b_i, c_ix+d_i)\mid f(x)(a_ix+b_i)^{l-m_i}-g(x)(c_ix+d_i)^{l-m_i}.
\end{array} \eqno(5.9)
$$
Assume that for some $i_0$ both polynomials on the r.h.s. of (5.9) are nonzero. Since deg$\Phi_l(a_{i_0}x+b_{i_0}, c_{i_0}x+d_{i_0})=l-1,$ the divisibilities (5.9) imply that
$$l-1\leq \theta+ m_{i_0},\ \ l-1 \leq \theta +l-m_{i_0}.$$
Adding these inequalities we get $2(l-1)\leq 2 \theta +l,$ hence $l \leq 2\theta +2,$ and $l \leq 2\theta +1,$ since $l$ is an odd prime.

To finish the proof we have to exclude the possibility that for every $i=1,2,\ldots, n$ at least one of the polynomials on the r.h.s. of (5.9) vanishes. Since $n\geq 2,$ there is $j\neq i, 1\leq j \leq n.$

Thus it is sufficient to prove that at most one of the polynomials
$$F_1=f(x)(c_ix+d_i)^{m_i}-g(x)(a_ix+b_i)^{m_i},$$
$$\ \ \ \ \ F_2=f(x)(a_ix+b_i)^{l-m_i}-g(x)(c_ix+d_i)^{l-m_i},$$
$$\ F_3=f(x)(c_jx+d_j)^{m_j}-g(x)(a_jx+b_j)^{m_j},$$
$$\ \ \ \ \ \ F_4=f(x)(a_jx+b_j)^{l-m_j}-g(x)(c_jx+d_j)^{l-m_j}$$
vanishes. Assume that at least two of these polynomials vanish. We consider several cases.

1) $F_1=F_2=0.$
(In the case $F_3=F_4=0$ we proceed similarly, replacing $i$ by $j.$)

From $f(x)(c_ix+d_i)^{m_i}=g(x)(a_ix+b_i)^{m_i}$ and $(f(x), g(x))=(a_ix+b_i, c_ix+d_i)=1$ it follows that
$$f(x)=\alpha (a_ix+b_i)^{m_i},\ \ g(x)=\alpha (c_ix+d_i)^{m_i}\ \ \mbox{for some}\ \alpha \in F^{*}. \eqno(5.10)$$
Analogously $f(x)(a_ix+b_i)^{l-m_i}=g(x)(c_ix+d_i)^{l-m_i}$ implies that
$$f(x)=\alpha^{\prime} (c_ix+d_i)^{l-m_i},\ \ g(x)=\alpha^{\prime} (a_ix+b_i)^{l-m_i}\ \ \mbox{for some}\ \alpha^{\prime} \in F^{*}. \eqno(5.11)$$
From (5.10) we get max(deg$f(x),$ deg$g(x)$)$=m_i$ and from (5.11) max(deg$f(x),$ deg$g(x)$)$=l-m_i.$ Hence $m_i=l-m_i$, so $l=2m_i,$ this is impossible, since $l$ is an odd prime.

2) $F_1=F_3=0.$ (In the case $F_2=F_4=0$ we proceed analogously).

Similarly as above we get
$$f(x)=\alpha (a_ix+b_i)^{m_i},\ \ g(x)=\alpha (c_ix+d_i)^{m_i},$$
$$f(x)=\alpha^{\prime} (a_jx+b_j)^{m_j},\ \ g(x)=\alpha^{\prime} (c_jx+d_j)^{m_j},$$
where $\alpha, \alpha^{\prime}\in F^{*}.$

Hence max(deg$f(x),$ deg$g(x))=m_i=m_j=:m.$ Therefore,
$$\frac{f(x)}{g(x)}=\Big(\frac{a_ix+b_i}{c_ix+d_i}\Big)^{m}=\Big(\frac{a_jx+b_j}{c_jx+d_j}\Big)^{m},$$
hence
$$\frac{a_ix+b_i}{c_ix+d_i}=\eta\cdot \frac{a_jx+b_j}{c_jx+d_j},$$
where $\eta^m=1, \eta\in F,$ thus $\eta\in W(F).$

It follows that
$$\left (\begin{matrix}
   a_i & b_i\\
   c_i & d_i
  \end{matrix}\right )=\lambda \left (\begin{matrix}
   \eta a_j & \eta b_j\\
   c_j & d_j
  \end{matrix}\right ),$$
where $\lambda\in F^*.$ This means that the matrices $\left (\begin{matrix}
   a_i & b_i\\
   c_i & d_i
  \end{matrix}\right )$ and $\left (\begin{matrix}
   a_j & b_j\\
   c_j & d_j
  \end{matrix}\right )$ are not essentially distinct. We get a contradiction, since $i\neq j.$

3) $F_1=F_4=0.$ (The case $F_2=F_4=0$ is quite analogous).

Similarly as above we get
$$f(x)=\alpha (a_ix+b_i)^{m_i},\ \ g(x)=\alpha (c_ix+d_i)^{m_i},$$
$$f(x)=\alpha^{\prime} (c_jx+d_j)^{l-m_j},\ \ g(x)=\alpha^{\prime} (a_jx+b_j)^{l-m_j},$$
where $\alpha, \alpha^{\prime}\in F^{*}.$

Hence max(deg$f(x),$ deg$g(x))=m_i=l-m_j=:m.$ Therefore,
$$\frac{f(x)}{g(x)}=\Big(\frac{a_ix+b_i}{c_ix+d_i}\Big)^{m}=\Big(\frac{c_jx+d_j}{a_jx+b_j}\Big)^{m},$$
hence
$$\frac{a_ix+b_i}{c_ix+d_i}=\eta\cdot \frac{c_jx+d_j}{a_jx+b_j},$$
where $\eta^m=1, \eta\in F,$ thus $\eta\in W(F).$

It follows that
$$\left (\begin{matrix}
   a_i & b_i\\
   c_i & d_i
  \end{matrix}\right )=\lambda \left (\begin{matrix}
   \eta c_j & \eta d_j\\
   a_j & b_j
  \end{matrix}\right )=\lambda \left (\begin{matrix}
   \eta & 0\\
   0 & 1
  \end{matrix}\right )\left (\begin{matrix}
   0 & 1\\
   1 & 0
  \end{matrix}\right ) \left (\begin{matrix}
   a_j \ & \  b_j\\
   c_j & d_j
  \end{matrix}\right ),$$
where $\lambda\in F^*.$ This means that the matrices $\left (\begin{matrix}
   a_i & b_i\\
   c_i & d_i
  \end{matrix}\right )$ and $\left (\begin{matrix}
   a_j & b_j\\
   c_j & d_j
  \end{matrix}\right )$ are not essentially distinct. We get a contradiction, since $i\neq j.$ \hfill$\Box$

\bigskip

\noindent {\bf Lemma 5.3}\ \ {\it Let ch$(F)=p> 0$ and $f,g\in F[x].$ Then:}

($i$) {\it If $f\notin F[x^p]$ and $f^r\in F[x^p],$ then $p|r.$}

($ii$) {\it If $(f,g)=1$ and $fg \in F[x^p],$ then $f,g\in F[x^p].$}

($iii$) {\it $F(x^p)\cap F[x]=F[x^p].$}

\noindent {\it Proof:} ($i$) By assumption $(f^r)^\prime =0$ and $f^\prime \neq 0.$ On the other hand, $(f^r)^\prime=rf^\prime f^{r-1}.$
Hence $r=0$ in $F,$ so $p|r.$

($ii$) We have $(fg)^\prime =0,$ hence $fg^\prime+f^\prime g=0.$ From $(f,g)=1$ it follows that $f|f^\prime$ and $g|g^\prime,$ then $f^\prime =g^\prime =0,$ that is $f,g \in F[x^p].$

($iii$) This formula is obvious. \hfill$\Box$

\bigskip

\noindent {\bf Lemma 5.4}\ \ {\it Assume that ch$(F)=p>0.$ If the polynomials $f,g$ in} (5.3) {\it belong to $F[x^p],$ then $\Psi \in F[x^p]$ and $p|r_i$ for every $i.$ Therefore} (5.3)
{\it implies an analogous formula with $x^p$ replaced by $x.$}

\noindent {\it Proof:} Let $f(x)=f_0(x^p), g(x)=g_0(x^p),$ where $f_0, g_0\in F[x].$ Then
$$\Phi_l(f(x),g(x))=\Phi_l(f_0(x^p), g_0(x^p))\in F[x^p].$$
By Lemma 5.3 ($ii$) and (5.3), the polynomials $\Psi(x)$ and $\Phi_l(a_ix+b_i, c_ix+d_i)^{r_i}$ belong to $F[x^p].$
 Thus $\Psi(x)=\Psi_0(x^p),$ where $\Psi_0 \in F[x].$

Since $\Phi_l(a_ix+b_i, c_ix+d_i)\notin F[x^p],$ then, by Lemma 5.3 ($i$), $p|r_i.$ So $r_i=pr_{i0}.$ We have
$\Phi_l(a_ix+b_i, c_ix+d_i)^p=\Phi_l((a_ix+b_i)^p, (c_ix+d_i)^p),$ because the polynomial $\Phi_l(x,y)$ has coefficients in $\mathbb{Z}/p.$

Obviously, $(a_ix+b_i)^p=a_{i0}x^p+b_{i0},$ and $(c_ix+d_i)^p=c_{i0}x^p+d_{i0},$ where $a_{i0},b_{i0}, c_{i0}, d_{i0} \in F. $
Therefore $\Phi_l(a_ix+b_i, c_ix+d_i)^p=\Phi_l((a_ix+b_i)^p=\Phi_l(a_{i0}x^p+b_{i0}, c_{i0}x^p+d_{i0}).$

Thus (5.3) can be written in the form of
$$\Phi_l(f_0(x^p), g_0(x^p))=\alpha \Psi_0(x^p)^l \prod_{i=1}^n\Phi_l(a_{i0}x^p+b_{i0}, c_{i0}x^p+d_{i0})^{r_{i0}}.$$
Replacing here $x^p$ by $x$ we get the formula analogous to (5.3).  \hfill$\Box$

\bigskip

\noindent {\bf Theorem 5.5}\ \ {\it Assume that ch$(F)=p$ and $f,g\in F[x^p].$ Then} (5.3) {\it can be written in the form of
$$\Phi_l(\widetilde{f}(x^{p^{r}}), \widetilde{g}(x^{p^{r}}))=\widetilde{\alpha} \widetilde{\Psi}(x^{p^r})^l\prod_{i=1}^n\Phi_l(\widetilde{a_{i}}x^{p^{r}}+\widetilde{b_{i}}, \widetilde{c_{i}}x^{p^{r}}+\widetilde{d_{i}})^{\widetilde{r}_{i}},$$
where $f(x)=\widetilde{f}(x^{p^{r}}), g(x)=\widetilde{g}(x^{p^{r}})$ and $\Psi(x)=\widetilde{\Psi}(x^{p^r})$ with  $\widetilde{f}^\prime \neq 0$ or $\widetilde{g}^\prime \neq 0,$ and $\widetilde{r}_i=r_i/p^{r} \in \mathbb{N}, \widetilde{\alpha}, \widetilde{a_{i}}, \widetilde{b_{i}}, \widetilde{c_{i}}, \widetilde{d_{i}}\in F.$}

\noindent {\it Proof:} If $f,g \in F[x^{p^{r}}],$ but at least one of these polynomials does not belong to $F[x^{p^{r+1}}],$ then applying Lemma 5.4 $r$ times
we get the formula analogous to (5.3) with the r.h.s. of the form $\Phi_l(\widetilde{f}, \widetilde{g}),$ where $f(x)=\widetilde{f}(x^{p^r}), g(x)=\widetilde{g}(x^{p^r}).$ Moreover, at least one of the polynomials $\widetilde{f}$ and $\widetilde{g}$ does not belong to $F[x^p],$ so  $\widetilde{f}^\prime \neq 0$ or $\widetilde{g}^\prime \neq 0.$  \hfill$\Box$

\bigskip

Let $\widetilde{f}, \widetilde{g}, \widetilde{\Psi}$ be as in Theorem 5.5, and  let $\widetilde{\theta}:=$max(deg$\widetilde{f},$ deg$\widetilde{g}$) and $\widetilde{\lambda}:=$deg$\widetilde{\Psi}.$ Then $\theta=p^r\cdot \widetilde{\theta}$ and $\lambda=p^r\cdot\widetilde{\lambda}.$ Note that $r_i=\widetilde{r}_i\cdot p^r.$

\bigskip

\noindent {\bf Theorem 5.6}\ \ {\it In the above notation, we have the inequalities}

(i) {\it If $f^{\prime}(x)\neq 0$ or $g^{\prime}(x)\neq 0,$ then}
$$n\leq \theta \leq \frac{(l-1)^2n-2l}{(l-1)^2-2l}.$$

(ii) {\it If $f^{\prime}(x)=g^{\prime}(x)= 0,$ then}
$$n\leq \widetilde{\theta} \leq \frac{(l-1)^2n-2l}{(l-1)^2-2l}.$$

(iii) {\it If $n\leq \frac{1}{2}(l^2-4l+1),$ then } deg$\Psi=0,$ {\it i.e. $\Psi\in F^*.$}

\noindent {\it Proof:}  (i)  Assume that $f^{\prime}(x)\neq 0$ or $g^{\prime}(x)\neq 0.$ Denote $\lambda:=$deg$\Psi, \theta:=$deg$f(x)\geq$deg$g(x).$ Then (5.3) implies
$$(l-1)\theta =l\lambda+(l-1)\sum_{i=1}^{n}r_i.\eqno(5.12)$$
Multiplying by $f(x)-g(x)$ both sides of (5.3) we get
$$
f^l - g^{l} = \alpha(f - g) \Psi^l \prod _{i=1}^{n}\Phi_l(a_ix+b_i,
c_ix+d_i)^{r_{i}}. \eqno(5.13)
$$

By the well known property of the differentiation, we have
$$\mbox{If}\ \ a,b\in F[x]\ \mbox{satisfy}\ a^r|b, r\geq 1,\ \mbox{then}\ a^{r-1}|b'.$$

Consequently, from (5.13) we get
$$F:= \Psi^{l-1}\prod_{i=1}^{n}\Phi_l(a_ix+b_i,c_ix+d_i)^{r_i-1}  \mid (f^l-g^l)^{\prime}=l(f'f^{l-1} - g'g^{l-1}). \eqno(5.14)
$$
By (5.13) ,
$$F\mid f^l-g^l. \eqno(5.15)$$
Hence $(F,f)=(F,g)=1,$ because $(f,g)=1.$

From
$$ g^{\prime}(f^{l}-g^{l})-g\cdot (f'f^{l-1} - g'g^{l-1})=f^{l-1}(fg^{\prime}-gf^{\prime})$$
and $(F,f)=1,$ by (5.14) and (5.15), we conclude that
$$F \mid fg^{\prime} - gf^{\prime}. \eqno(5.16)
$$

Since $f^{\prime}(x)\neq 0$ or $g^{\prime}(x)\neq 0,$
 we get $fg^\prime-gf^\prime \neq 0.$
Therefore from (5.14) and (5.16) it follows that
$$\mbox{deg}F=(l-1)\lambda +(l-1)\sum_{i=1}^{n}(r_i-1)\leq \mbox{deg}(fg^\prime-gf^\prime)\leq 2\theta -2. \eqno(5.17)$$
Namely, it is an easy exercise to prove that for any polynomials $f,g\in F[x]$ satisfying $\theta=$deg$f\geq$ deg$g$ and $fg^\prime-gf^\prime \neq 0$
we have deg$(fg^\prime-gf^\prime)\leq 2 \theta -2.$ It is sufficient to consider the leading terms of $f$ and $g.$

Thus we have proved the two formulas (5.12)
 and (5.17) relating $l,\lambda$ and $\theta.$ From (5.12) it follows that $l-1|\lambda,$ so $\lambda=(l-1)\lambda_1,$ where $\lambda_1\geq 0.$

Dividing (5.12) and (5.17) by $l-1$ we get
$$\theta =l\lambda_1+\sum_{i=1}^nr_i, \eqno(5.18)$$
$$(l-1)\lambda_1+\sum_{i=1}^nr_i-n\leq \frac{2}{l-1}(\theta-1).\eqno(5.19)$$
Since $r_i \geq 1$ and $1\leq i \leq n,$ then $\sum_{i=1}^nr_i\geq n.$ Consequently, (5.18) and (5.19) imply
$$\theta \geq l\lambda_1+n, \eqno(5.20)$$
$$(l-1)\lambda_1\leq \frac{2}{l-1}(\theta-1). \eqno(5.21)$$
From (5.20) we get $\theta \geq n,$ which gives the first inequality in $(i).$

By (5.18), (5.19) and (5.21), we have
$$\theta =\sum_{i=1}^nr_i+(l-1)\lambda_1+\lambda_1 \leq n+\frac{2}{l-1}(\theta-1)+\frac{2}{(l-1)^2}(\theta-1)=n+\frac{2l}{(l-1)^2}(\theta-1).$$
Hence
$$\theta \leq \frac{(l-1)^2n-2l}{(l-1)^2-2l}.$$
This gives the second inequality in $(i).$

(ii) Assume that  $f^{\prime}(x)= 0$ and $g^{\prime}(x)= 0.$ Clearly we must have ch$(F)=p> 0$ and  $f, g\in F[x^p].$
By Theorem 5.5, we have
$$\Phi_l(\widetilde{f}(X), \widetilde{g}(X)=\widetilde{\alpha} \widetilde{\Psi}^l\prod_{i=1}^n\Phi_l(\widetilde{a_{i}}X+\widetilde{b_{i}}, \widetilde{c_{i}}X+\widetilde{d_{i}})^{\widetilde{r}_{i}},$$
where $f(x)=\widetilde{f}(X), g(x)=\widetilde{g}(X), X=x^{p^{r}}$ with  $\widetilde{f}^\prime(X) \neq 0$ or $\widetilde{g}^\prime(X) \neq 0,$

Let $\widetilde{\lambda}=(l-1)\widetilde{\lambda}_1.$ Then similarly we have
$$\widetilde{\theta} =l\widetilde{\lambda}_1+\sum_{i=1}^n\widetilde{r}_i, \eqno(5.22)$$
$$(l-1)\widetilde{\lambda}_1+\sum_{i=1}^n\widetilde{r}_i-n\leq \frac{2}{l-1}(\widetilde{\theta}-1).\eqno(5.23)$$
Since  $\sum_{i=1}^n\widetilde{r}_i\geq n,$ consequently, (5.22) and (5.23) imply
$$\widetilde{\theta} \geq l \widetilde{\lambda}_1+n, \eqno(5.24)$$
$$(l-1)\widetilde{\lambda}_1\leq \frac{2}{l-1}(\widetilde{\theta}-1). \eqno(5.25)$$
From (5.24) we get $\widetilde{\theta} \geq n,$ which gives the first inequality in $(ii).$

By (5.22), (5.23) and (5.25), we have
$$\widetilde{\theta} \leq \frac{(l-1)^2n-2l}{(l-1)^2-2l}.$$
This gives the second inequality in $(ii).$

(iii) If $f^{\prime}(x)\neq 0$ or $g^{\prime}(x)\neq 0,$ from (5.21) and $(i)$ we obtain
$$\lambda_1 \leq \frac{2}{(l-1)^2}(\theta -1)\leq \frac{2(n-1)}{l^2-4l+1}.$$
It follows that $\lambda_1< 1$ if $n-1< \frac{1}{2}(l^2-4l+1).$ Since $l^2-4l+1$ is an integer, the last inequality is equivalent to $n\leq \frac{1}{2}(l^2-4l+1).$
This proves  $\lambda_1=0,$ so deg$\Psi=\lambda=(l-1)\lambda_1=0.$

If $f^{\prime}(x)=g^{\prime}(x)= 0,$
from (5.25) and $(ii)$ we obtain
$$\widetilde{\lambda}_1 \leq \frac{2}{(l-1)^2}(\widetilde{\theta} -1)\leq \frac{2(n-1)}{l^2-4l+1}.$$
Similarly $n\leq\frac{1}{2}(l^2-4l+1)$ implies  that $\widetilde{\lambda}_1=0,$ that is deg$\Psi=(l-1)\widetilde{\lambda}_1 p^r=0.$
 \hfill$\Box$

\noindent {\it Remarks} 5.7 \ i) The argument above is analogous to the proof of the $abc$-conjecture for polynomials.

{\bf Theorem abc} (W.W. Stothers). {\it Let $a,b,c\in F[x],$  where} char$F=0$ {\it and not all polynomials $a,b,c$ are constant. For a nonzero polynomial $h\in F[x]$ denote by} rad$(h)$ {\it the number of distinct roots of $h$ in the algebraic closure of $F.$ Assume that $a,b,c$ are relatively prime and $a+b=c.$ Then}
$$\mbox{max}(\mbox{deg}a, \mbox{deg}b, \mbox{deg}c)\leq \mbox{rad}(abc)-1.$$

We can apply the theorem $abc$ as follows. In the notation of (5.13) put $a:=f^l, b:=-g^{l},$ and $c:=$ the r.h.s. of (5.13).

Then max(deg$a,$ deg$b,$ deg$c$)$=$deg$(f^l)=l\theta,$ rad$(a)\leq $deg$f=\theta,$ rad$(b)\leq$deg$g\leq \theta,$ and
$$\mbox{rad}(c)\leq \mbox{deg}(f-g)+\mbox{deg}\Psi+\sum_{i=1}^{n}\mbox{deg}\Phi_l(f_i,g_i)\leq \theta+\lambda+n(l-1).$$
Consequently, the theorem $abc$ gives
$$l\theta \leq 3 \theta+(l-1)\lambda_1+(l-1)n-1.$$
Considering all terms of this inequality modulo 2, we see that the last term $-1$ can be replaced by $-2.$

Hence
$$\theta \Big(1-\frac{2}{l-1}\Big)\leq \lambda_1+n-\frac{2}{l-1}.$$
Now, applying the estimate $\lambda_1 \leq \frac{2}{(l-1)^2}(\theta-1)$ following from (5.21), we get
$$\theta\leq \frac{(l-1)^2n-2l}{(l-1)^2-2l}.$$
 Thus we obtain the second inequality in Theorem 5.6 $(i).$

ii) When $n=1,$ (5.3) is trivial. In fact, we can prove the following statement:

 {\it We still assume that $f^\prime \neq 0$ or $g^\prime \neq 0.$ In the case $n=1$ the formula} (5.3) {\it takes the form $\Phi_l(x)=\Phi_l(x).$}

{\it Proof.} By Theorem 5.6 $(i),$ from $n=1$ it follows that $\theta =1,$ that is deg$f=1\geq$deg$g.$ Hence $f(x)=ax+b, g(x)=cx+d,$ where $a=1,$ since we always assume that $f$ is monic.

Therefore (5.3) takes the form
$$\Phi_l(f,g)=\alpha \Phi_l(x)^{r_1}\ \ \mbox{for some}\ \alpha \in F^*,\eqno(5.26)$$
since deg$\Psi=0,$ by Theorem 5.6 $(iii).$

Comparing degrees of both sides of (5.26) we get $r_1=1.$ From (5.26) it follows that the polynomials $\Phi_l(f,g)=\Phi_l(ax+b,cx+d)$ and $\alpha \Phi_l(x)$ are not relatively prime.

Then, by Theorem 2.5, the corresponding matrices
$$\left (\begin{matrix}
   1 & b\\
   c & d
  \end{matrix}\right )\ \ \mbox{and}\ \ \left (\begin{matrix}
   1 & 0\\
   0 & 1
  \end{matrix}\right )$$
are not essentially distinct. Therefore from $a=1$ it follows that $b=c=0$ and $d=\mu$ is a root of unity. Hence
$$\Phi_l(f,g)=\Phi_l(x,\mu)=\alpha\Phi_l(x).$$
Comparing the leading terms we get $\alpha =1,$ then the coefficients by $x^{l-2}$ in both polynomials are $1$ and $\mu.$ Hence $\mu=1,$ so $f(x)=x, g(x)=1,$ and (5.26) takes the form $\Phi_l(x)=\Phi_l(x).$  \hfill$\Box$

\bigskip

\noindent {\bf Theorem 5.8}\ {\it In the above notation, assume that  that  $2\leq n\leq \frac{1}{2}(l^2-4l+1).$  Then we have}
$$l\leq 2n+1.$$

\noindent {\it Proof:}  Assume that $f^{\prime}(x)\neq 0$ or $g^{\prime}(x)\neq 0.$ Then in the case $\theta=n$ from Theorem 5.2 it follows that $l\leq 2\theta +1=2n+1.$

If $\theta> n, $  from Theorem 5.6 $(iii)$ we get $\lambda_1=0,$ then (5.18) and (5.19) give $1\leq \theta-n \leq \frac{2}{l-1}(\theta-1).$
Hence
$$l\leq 1+2\cdot \frac{\theta-1}{\theta-n}=3+2\cdot \frac{n-1}{\theta-n}\leq 3+2(n-1)=2n+1.$$

 Assume that  $f^{\prime}(x)=g^{\prime}(x)= 0.$ Then we have $\widetilde{f}^{\prime}(x)\neq 0$ or $\widetilde{g}^{\prime}(x)\neq 0.$ In the case $\widetilde{\theta}=n$ from Theorem 5.2 it follows that $l\leq 2\widetilde{\theta} +1=2n+1.$

If $\widetilde{\theta}> n, $  from Theorem 5.6 $(iii)$ we get $\widetilde{\lambda}_1=0,$ then (5.22) and (5.23) give $1\leq \widetilde{\theta}-n \leq \frac{2}{l-1}(\widetilde{\theta}-1).$
Hence $l\leq 2n+1.$  \hfill$\Box$

\bigskip

\noindent {\bf Corollary 5.9}\  {\it Assume that $l\geq 5$ is a prime number and  $F$ is a field
such that $\Phi_l(x)$ is irreducible in $F[x]$. Let $n$ be an integer satisfying
$$2\leq n\leq \frac{l-3}{2},$$ and let
$\gamma_1, \gamma_2,\ldots, \gamma_n \in GG_l(F(x))$ be essentially distinct. Then
$$\prod_{i=1}^{n}\gamma_i^{l_i} \notin G_l(F(x)),$$
where $1\leq l_i \leq l-1, i=1,2, \ldots,n.$}

\noindent {\it Proof:} It follows from  Theorem 5.5 and 5.8.  \hfill$\Box$

\bigskip

\noindent {\bf Corollary 5.10}\ {\it Assume that $l\geq 5$ is a prime number and  $F$ is a field
such that $\Phi_l(x)$ is irreducible in $F[x]$. Let $n$ be an integer satisfying
$$2\leq n\leq \frac{l-3}{2}.$$
If $H$ is a cyclotomic subgroup of $\mathfrak{G}_{l}(n; F)$, then it is a cyclic group of order $l.$}  \hfill$\Box$

\bigskip

The following result gives the  relations between $n$ and $\theta$ (or $\widetilde{\theta}$).

\bigskip

\noindent {\bf Theorem 5.11}\ \ {\it Assume that $n\leq \frac{1}{2}(l^2-4l+1).$ Then we have:}

(i) {\it If $f^{\prime}(x)\neq 0$ or $g^{\prime}(x)\neq 0,$ then }
 $ \theta \leq 2n-1.$

 (ii) {\it If $f^{\prime}(x)=g^{\prime}(x)= 0,$ then} $\widetilde{\theta} \leq 2n-1.$

\noindent {\it Proof:}\ If $\theta> 2n-1,$ then from Theorem 5.6 (iii), (5.18)and  (5.19) we have
 $$l\leq 3+2\cdot \frac{n-1}{\theta-n}< 3+2\cdot \frac{n-1}{(2n-1)-n}=5,$$
  which contradicts the assumption that $l\geq 5.$ Hence $\theta \leq 2n-1.$ Then proof of $\widetilde{\theta} \leq 2n-1$ is similar.  \hfill$\Box$

\bigskip
\noindent {\it Remarks} 5.12 \
a) More precisely, in the case $f^{\prime}(x)\neq 0$ or $g^{\prime}(x)\neq 0,$ from the proof of Theorem 5.11 it follows that

(i)  if $\theta=n,$ then $l\leq 2n+1;$

(ii)  if $n< \theta \leq 2n-1,$ then $l\leq 3+2\cdot \frac{n-1}{\theta-n}.$

In particular, we have
$$\mbox{If}\ \theta=n+1, \ \mbox{then}\ l\leq 2n+1,$$
$$\mbox{If}\ \theta=n+2, \ \mbox{then}\ l\leq n+2,$$
$$\mbox{If}\ \theta=2n-1, \ \mbox{then}\ l\leq 5.\ \ \ \ \ $$

b) From Theorem 5.2, we get the relation between $l$ and $\theta,$ i.e.,
$l\leq 2\theta+1.$ Furthermore, if ch$(F)=0$ and $\theta > n,$ then from (5.19) we have $l\leq 2\theta-1.$
As was suggested to me by Browkin, the last inequality is actually a necessary condition for the polynomial $\Phi_l(f,g)$  having a multiple root. In fact, we can prove the following statement:

\bigskip

 {\it Assume that $l\geq 5$ is a prime and $\Phi_l(x)$ is irreducible in $F[x].$ If ch$(F)=0$ and $\Phi_l(f,g)$ has a multiple root, where gcd$(f,g)=1,$ then
$l\leq 2\theta-1.$

 In particular, if $l=5, \theta=2$ and $l=7, \theta=2$ or $3,$ then $\Phi_5(f,g)$ and $\Phi_7(f,g)$ have no a multiple root, respectively.}

\noindent {\it Proof:} Assume that $\alpha$ is a multiple root of $\Phi_l(f,g),$ then it must be a multiple root of $f(x)-\zeta g(x),$ where $\zeta=\zeta _l.$ Then
$$f(\alpha)-\zeta g(\alpha)=0,\ \ \ f^{\prime}(\alpha)-\zeta g^{\prime}(\alpha)=0,$$
so
$$f(\alpha)g^{\prime}(\alpha)-f^{\prime}(\alpha)g(\alpha)=0.$$

It follows that $\alpha$ is a root of the polynomial $t(x):=f(x)g^{\prime}(x)-f^{\prime}(x)g(x).$ From
$(f,g)=1$ and ch$(F)=0$ it follows that $t(x)$ is a non-zero polynomial of degree at most $2\theta-1.$

From $f(\alpha)-\zeta g(\alpha)=0$ we conclude that $F\subseteq F(\zeta)\subseteq F(\alpha).$ Since $[F(\alpha):F]$ is the degree of the minimal polynomial for $\alpha$ over $F,$ and $[F(\alpha):F]$ is divisible by $[F(\zeta):F]=l-1$ we conclude that $l-1\leq 2\theta-1,$ i.e., $l\leq 2\theta,$ so $l\leq 2\theta-1$ since $l$ is odd, as claimed.  \hfill$\Box$

\bigskip

Now, we turn to the case of $n=1.$
Let $l,p$ be two different prime numbers. Define
$$\mathfrak{Z}(l,p):=\{t \mid 2\leq t \leq l-2, \ t\equiv p^{2m} \mbox{or} -p^{2m} (\mbox{mod} l)\ \mbox{for some}\ m\in \mathbb{N} \}.$$

\bigskip

\noindent {\bf Lemma 5.13}\ \ {\it Assume that $l\geq 5$ is a prime number and  $F$ is a field
such that $\Phi_l(x)$ is irreducible in $F[x]$. Let $\gamma\in GG_l(F(x)).$

i) If ch$(F)=0, $ then none of the elements $\gamma^{t}, 2\leq t\leq
l-2,$ is  cyclotomic. So, the only cyclotomic elements contained in $\langle \gamma\rangle$ are $\gamma, \gamma^{-1}.$
Hence, $\langle \gamma\rangle$ is not a cyclotomic subgroup.

ii) If ch$(F)=p\neq0,$ then
$$1\neq\gamma^t\in G_l(F(x)) \Longleftrightarrow t\in \{1, l-1\} \cup \mathfrak{Z}(l,p).$$
So $\langle \gamma\rangle$ contains exactly $2+|\mathfrak{Z}(l,p)|$ nontrivial cyclotomic elements.
 }

\noindent {\it Proof:}\ Clearly, it suffices to
consider $\gamma=  c_l(x).$ Let $t$ be a temporarily
fixed integer satisfying $2\leq t \leq l-2.$ If $ \gamma^{t} $ is
cyclotomic, then there exist nontrivial polynomials $f_t,g_t\in
F[x]$ such that
$$ \gamma^{t} =c_l \left (\frac{f_t}{g_t}\right),$$
with $f_t$  monic. By Theorem 5.1 ($i$), we have the
equality:
$$
\Phi_l(f_t,g_t) = \alpha_t \Psi_t^l \Phi_l(x)^{r_t}. \eqno(5.27)
$$

Let $\theta_t:=$max(deg$f_t,$ deg$g_t$) and $\lambda_t:=$deg$\Psi_t.$

i) Assume that ch$(F)=0.$ Then $f^{\prime}(x)\neq 0$ or $g^{\prime}(x)\neq 0.$ From  Theorem 5.6 (i) we have
 $\theta_t=1,$ hence $ \lambda_t=0$ and $r_t=1.$

Now, let
$$f_t(x)=a_tx+b_t, \ \ g_t(x)=c_tx+d_t.$$
Then (5.27) becomes
$$\Phi_l(a_tx+b_t, c_tx+d_t)=\alpha_t \Phi_l(x).$$

Let $x=\zeta.$ Then there exists an $i$ satisfying $1\leq i\leq l-1$ such that
$$\frac{a_t\zeta+b_t}{c_t\zeta+d_t}=\zeta^i,$$
so
$$c_t\zeta^{i+1}+d_t\zeta^{i}-a_t\zeta-b_t=0.$$
Easy computations show that the possible cases are only either
$a_t=d_t\neq0, b_t=c_t=0$ or $b_t=c_t\neq 0, a_t=d_t=0. $  So
$$f_t(x)=a_tx, \ \ g_t(x)=a_t $$
or
$$  f_t(x)=c_t, \ \  g_t(x)=c_tx.$$

 Hence, if
$f_t(x)=a_tx, g_t(x)=a_t$ we get
$$c_l(x)^{t}=\beta=c_l\left(\frac{f_t}{g_t}\right)=c_l(x),$$
which implies $c_l(x)=1,$ a contradiction; if
$f_t(x)=c_t, g_t(x)=c_tx,$ we get
$$c_l(x)^{t}=c_l(x^{-1})=c_l(x)^{-1},$$
so $c_l(x)^{t+1}=1,$ therefore $c_l(x)=1$ since
$2\leq t \leq l-2,$ also a contradiction.

In summery, the equality (5.27) does not hold. So  none of
$\gamma^{t}, 2\leq t \leq l-2,$ is cyclotomic.

ii) Assume that ch$(F)=p>0.$ If there exists  some $t$ satisfying $ 2\leq t \leq l-2$ such
that $f_t^{\prime}\neq0$ or $ g_t^{\prime}\neq 0,$ then discussions
similar as i) show that $c_l(x)^{t}$ is not cyclotomic.
Hence, if $\{x, \Phi_l(x)\}^{t}$ is cyclotomic for some  $2\leq t \leq l-2,$ we must have
$f_t^\prime=0$ and $ g_t^\prime=0.$

Similarly as in (i), we have
$$\Phi_l(a_tx^{p^{m_t}}+b_t, c_tx^{p^{m_t}}+d_t)=\alpha_t \Phi_l(x^{p^{m_t}}).$$

Let $x^{p^{m_t}}=\zeta.$ Then we get
$\frac{a_t\zeta+b_t}{c_t\zeta+d_t}=\zeta^i$ for some $i$ satisfying
$1\leq i \leq l-1.$ A computation leads to either $a_t=d_t,
b_t=c_t=0$ or $a_t=d_t=0, b_t=c_t.$ So  we have either
$$f_t(x)=a_tx^{p^{m_t}},\ \  g_t(x)=a_t$$
or
$$f_t(x)=a_t,\ \  g_t(x)=a_tx^{p^{m_t}}.$$

Therefore, if $f_t(x)=a_tx^{p^{m_t}},  g_t(x)=a_t,$ we have
$$c_l(x)^{t}=\beta=c_l(x^{p^{m_t}})=c_l(x)^{p^{2m_t}}.$$
Hence
$l\mid p^{2m_t}-t,$ that is  $t\in \mathfrak{Z}(l,p);$
 if $f_t(x)=a_t,  g_t(x)=a_tx^{p^{m_t}},$ then
$l\mid p^{2m_t}+t,$ also $t\in \mathfrak{Z}(l,p).$
Hence, for $2\leq t\leq l-2,$ if $c_l(x)^{t}$ is cyclotomic, then we have $t\in \mathfrak{Z}(l,p).$

On the other hand, if $t\in \mathfrak{Z}(l,p),$ then
we have either
$$t=p^{2m_t}+lm^\prime,\ \ \mbox{for some integer} \ m^\prime,$$
or
$$t=-p^{2m_t}+lm^{\prime\prime}, \ \ \mbox{for some integer} \ m^{\prime\prime}.$$
So we have either
$$c_l(x)^{t}=c_l(x)^{p^{2m_t}+lm^\prime}=c_l(x^{p^{m_t}})$$
or
$$c_l(t)^{t}=c_l(x)^{-p^{2m_t}+lm^{\prime\prime}}=c_l(x)^{-p^{2m_t}}=c_l(x^{-p^{m_t}}).$$
This implies that if $t\in \mathfrak{Z}(l,p),$ then we have $c_l(x)^{t}\in G_l(F(x)).$

Note that $c_l(x), c_l(x)^{-1}\in G_l(F(x)).$ Then we get the lemma.
 \hfill $\Box$

\bigskip

\noindent {\bf Lemma 5.14}\  {\it  The following statements are equivalent.

i) $|\mathfrak{Z}(l,p)|= l-3.$

ii)  $ l\equiv 3\, (\mbox{mod}\, 4)$ and
$p$ is a primitive root of $l$.

}

\noindent {\it Proof:}  Clearly, if $p$ is not a primitive root of $l,$ then the order of $p^2$(mod $l$) is less than $\frac{l-3}{2}.$
So $|\mathfrak{Z}(l,p)|<l-3.$

 When $p$ is a primitive root of $l,$ the set of all quadratic residues  (mod $l$) is
$$1, p^2, p^4, \ldots, p^{2(\frac{l-3}{2})}.$$

Consider the map: $p^{2m}\mapsto -p^{2m}.$ This is a bijection. If $l\equiv 3\, (\mbox{mod}\, 4),$ then we have
$$\left( \frac{-p^{2m}}{l}\right)=\left( \frac{-1}{l}\right)=(-1)^{\frac{l-1}{2}}=-1,$$
where $\left(\frac{\cdot}{l}\right)$ is the Legendre symbol (mod$l$).
Hence, if $t\equiv -p^{2m} (\mbox{mod}\, l),$ then $t$ is a quadratic non-residue (mod $l$).
So $| \mathfrak{Z}(l,p)|= l-3.$

Conversely, if $l\equiv 1\, (\mbox{mod}\, 4),$ then
$$\left( \frac{-p^{2m}}{l}\right)=\left( \frac{-1}{l}\right)=1.$$
This implies that the integers in $\mathfrak{Z}(l,p)$ are all quadratic residues (mod $l$). But the number of quadratic residues is $\frac{l-1}{2}.$ So
$$| \mathfrak{Z}(l,p)| \leq \frac{l-1}{2}< l-3,$$
a contradiction. Hence $l\equiv 3\, (\mbox{mod}\, 4).$  \hfill$\Box$

\bigskip

\noindent {\bf Corollary 5.15}\  {\it Assume that $l\geq 5$ is a prime number and  $F$ is a field with ch$(F)=p$
such that $\Phi_l(x)$ is irreducible in $F[x]$.
For any $\gamma \in GG_l(F(x)),$ the subgroup of $K_2(F(x))$ generated by
$\gamma$ is cyclotomic if and only if $l\equiv 3 \,(\mbox{mod}\, 4)$ and $p$ is a primitive root of $l,$ i.e.,
$$\langle \gamma \rangle \subset G_l(F(x)),\ \forall \ \gamma \in GG_l(F(x))\Longleftrightarrow  \ l\equiv 3\, (\mbox{mod}\, 4) \ \mbox{and}
\ p \mbox{ is a primitive root of} \ l.$$
 }
\noindent {\it Proof:}   \ Clearly, we have
$$\{\gamma, \gamma^{-1}\}\cup \{\gamma^{t} \mid t\in \mathfrak{Z}(l,p) \}\subseteq \langle \gamma \rangle,$$
which implies that
$$2+|\mathfrak{Z}(l,p)|=| \{\gamma, \gamma^{-1}\}\cup \{\gamma^{t} \mid t\in \mathfrak{Z}(l,p) \}| < | \langle \gamma \rangle | =l.$$

\bigskip

  If $l\equiv 3\, (\mbox{mod}\, 4),$ then from Lemma 5.14, we have $2+|\mathfrak{Z}(l,p)|=l-1,$ so from Lemma 5.13 ii), we get
$$\langle \gamma \rangle=\{1, \gamma, \gamma^{-1}\}\cup \{\gamma^{t} \mid t\in \mathfrak{Z}(l,p) \} \subseteq G_l(F(x)).$$

Conversely, from Lemma 5.13 ii), we have
$$\langle \gamma \rangle \subseteq \{1, \gamma, \gamma^{-1}\}\cup \{\gamma^{t} \mid t\in \mathfrak{Z}(l,p) \}\subseteq \langle \gamma \rangle.$$
So $l=3 +|\mathfrak{Z}(l,p)|,$ that is, $|\mathfrak{Z}(l,p)|=l-3.$ From Lemma 5.14, we have $l\equiv 3\, (\mbox{mod}\, 4).$
$\Box$

\bigskip

\noindent {\bf Example 5.16}\ \ It is easy to show that $\Phi_7(x)$ is irreducible in $\mathbb{F}_3[x]$ and $3$ is a primitive root of $7$.  \hfill $\Box$

\bigskip

Now we arrive at the main result of this section as follows.

\bigskip

\noindent {\bf Theorem 5.17}\  {\it Assume that $l\geq 5$ is a prime number and  $F$ is a field
such that $\Phi_l(x)$ is irreducible in $F[x]$. Let $n$ be an integer satisfying
$$ n\leq \frac{l-3}{2}.$$

 i) If ch$(F)=0, $ then $c(\mathfrak{G}_{l}(n; F))=2n,$  and so $cs(\mathfrak{G}_{l}(n; F))=0.$

 ii) If ch$(F)=p\neq0,$ then $c(\mathfrak{G}_{l}(n; F))=n(2+|\mathfrak{Z}(l,p)|).$

iii) If ch$(F)=p\neq0,$ then we have
  $$cs(\mathfrak{G}_{l}(n; F))> 0 \Longleftrightarrow  \ l\equiv 3\, (\mbox{mod}\, 4) \ \mbox{and}
\ p \mbox{ is a primitive root of} \ l.$$
   In this case,
  $cs(\mathfrak{G}_{l}(n; F))=n,$ i.e.,  $\mathfrak{G}_{l}(n; F)$ contains exactly $n$ nontrivial  cyclotomic subgroups.

iv) Every nontrivial cyclotomic subgroup of $\mathfrak{G}_{l}(n; F)$ is a cyclic subgroup of order $l,$ i.e.,
every nontrivial cyclotomic subgroup has the form $\mathfrak{G}_{l}(1; F).$}

\bigskip

\noindent {\it Proof:} \ i) It follows from Corollary 5.9 and Lemma 5.13 i).

ii) It follows from  Corollary 5.9 and Lemma 5.13 ii).

iii) It follows from Corollary 5.15.

iv) It follows from iii) and Corollary 5.10.
 \hfill$\Box$

\bigskip

\noindent {\bf Corollary 5.18}\  {\it Assume that $l\geq 5$ is a prime number with $l\equiv 3\, (\mbox{mod}\, 4)$ and  $F$ is a field with ch$(F)=p$
such that $\Phi_l(x)$ is irreducible in $F[x]$. If $p$  is a primitive root of $l,$ then $\mathfrak{G}_{l}(1; F)$ is a cyclotomic subgroup.} \hfill$\Box$

\bigskip

\noindent {\it Remark} 5.19\ \  From Theorem 5.17, we conclude immediately that  $G_l(F(x))$ is not a group, as is conjectured by Browkin in [1].

\bigskip

\noindent {\bf Corollary 5.20}\ \ {\it Assume that $l\geq 5$ is a prime number.   If $n$ is a positive integer satisfying
$$n\leq \frac{l-3}{2},$$
 then $c(\mathfrak{G}_l(n;\mathbb{Q}))=2n,$  so $cs(\mathfrak{G}_l(n;\mathbb{Q}))=0$.}  \hfill$\Box$

\bigskip

\noindent {\bf Corollary 5.21}\  {\it Assume that $l$ is a prime number, $F$ is a field
with ch$(F)\neq l$ and $\Phi_l(x)$ is irreducible in $F[x]$.

i) If ch$(F)=0 $ and  $l\geq 5$ (resp. $l\geq 7$ or $l\geq 11$), then $c(\mathfrak{G}_l(1;F))=2$ (resp. $c(\mathfrak{G}_l(2;F))=4$ or $c(\mathfrak{G}_l(3;F))=6$ and $c(\mathfrak{G}_l(4;F))=8$).

ii) If ch$(F)=p\neq0 $ and  $l\geq 5$ (resp. $l\geq 7$ or $l\geq 11$), then $c(\mathfrak{G}_l(1;F))=2+| \mathfrak{Z}(l,p)|$ (resp. $c(\mathfrak{G}_l(2;F))=2(2+|\mathfrak{Z}(l,p)|)$ or $c(\mathfrak{G}_l(3;F))=3(2+| \mathfrak{Z}(l,p)|)$ and $c(\mathfrak{G}_l(4;F))$
$ =4(2+|\mathfrak{Z}(l,p)|)).$
} \hfill$\Box$

\bigskip

\noindent {\bf Corollary 5.22}\  {\it Assume that $l\geq 5$ is a prime number, $F$ is a field
with ch$(F)\neq l$ and $\Phi_l(x)$ is irreducible in $F[x]$.  If $n$ is a positive integer satisfying
$$n\leq \frac{l-3}{2},$$

 i) If ch$(F)=0,$ then  $c(\mathfrak{S}(n;F))=c(\mathfrak{T}_l(n;F)=2n.$

ii) If ch$(F)=p\neq0, $ then
$$c(\mathfrak{S}_{l}(n; F))=c(\mathfrak{T}_l(n;F))=n(2+ |\mathfrak{Z}(l,p)|).$$
In particular, when $p$ is a primitive root of $l$ and $l\equiv 3\, (\mbox{mod}\, 4),$ we have
$$cs(\mathfrak{S}_{l}(n; F))=cs(\mathfrak{T}_l(n;F))=n.$$}
\hfill$\Box$

\bigskip

\noindent {\it Remark} 5.23 \ The equality (5.3) is actually a diophantine equation about $X,Y,Z$ over the polynomial ring $F[x],$ i.e., it can be rewritten as
$$\frac{X^l-Y^l}{X-Y}=\alpha \prod_{1}^{n}\Phi_l(a_ix+b_i,
c_ix+d_i)^{e_{i}}\cdot Z^l,$$
where $1\leq e_i\leq l-1$ and $a_id_i-b_ic_i=1$ with $ 1\leq i\leq n.$
If  $l\geq 5$ is a prime number and  $\Phi_l(x)$ is irreducible in $F[x]$, then from the proof of Theorem 5.17 we know that
the above diophantine equation  has no solution in $F[x]$ if $ n\leq \frac{l-3}{2}.$

\bigskip

Let $\mathbb{F}_q$ be a finite field of $q$ elements, where $q$ is a power of the prime $p>2,$ and for an integer $m>0,$ define
$$GG_l(\mathbb{F}_q(x))^{m}:=\{c^m: c\in GG_l(\mathbb{F}_q(x))\}.$$
\bigskip

\noindent {\bf Corollary 5.24}\ {\it Assume that $l\geq 5$ is a prime with $l\equiv 3\, (\mbox{mod}\, 4)$ and $l\neq p,$  that  $\Phi_l(x)$ is irreducible in $\mathbb{F}_p[x],$ and that $l,p$ satisfy the relation
$$n:=p(p+1)\leq l-3.$$

If $p$ is a primitive root of $l,$  then the set cyclotomic elements
$ \bigcup_{m=0}^{l-1} GG_l(\mathbb{F}_p(x))^{m}$ contains at least
 $n$  distinct nontrivial  cyclotomic subgroups, i.e., there are $n$ essentially distinct   elements $c_l\Big(\frac{a_ix+b_i}{c_ix+d_i}\Big), 1\leq i \leq n,$ so that
 $$G_l(\mathbb{F}_p(x))\supseteq \bigcup_{m=1}^{l-1} GG_l(\mathbb{F}_p(x))^{m}\supseteq\bigcup_{i=1}^{n}\Big\langle c_l\Big(\frac{a_ix+b_i}{c_ix+d_i}\Big) \Big\rangle.$$

}

 \noindent {\it Proof:} \ At first, since $l\equiv 3\, (\mbox{mod}\, 4)$ and   $p$ is a primitive root of $l,$ we have
 $$\bigcup_{m=1}^{l-1} GG_l(\mathbb{F}_p(x))^{m}=\bigcup_{t\in \{1,\, l-1\}\cup \mathfrak{Z}(l,p)}  GG_l(\mathbb{F}_p(x))^t \subseteq  G_l(\mathbb{F}_p(x)).$$

It is well known that $|PGL(2, \mathbb{F}_p)|=p(p^2-1).$ Hence from Lemma 4.3, we have
$$|GG_l(\mathbb{F}_p(x))|=|PGL(2, \mathbb{F}_p)|=p(p^2-1).$$

According to the definition,   if $A=\left (\begin{matrix}
   a & b\\
   c & d
  \end{matrix}\right )\in GL(2, \mathbb{F}_p)$ then all matrices which are not essentially distinct from $A$ are
  $$\alpha \left (\begin{matrix}
   \mu a & \mu b\\
   c & d
  \end{matrix}\right )\ \  \mbox{and}\ \ \alpha \left (\begin{matrix}
   \mu c & \mu d\\
   a & b
  \end{matrix}\right ), \ \mbox{for all}\ \alpha, \mu\in \mathbb{F}_p^*.\eqno(5.28)$$

Since $p>2,$ it is easy to show that the matrices of (5.28) are different from each other,
 so the number of elements in each class of non-essentially distinct elements is $2(p-1)^2.$
Therefore the number of classes of essentially distinct elements is
$$\frac{|GL(2,\mathbb{F}_p)|}{2(p-1)^2}=\frac{(p^2-1)(p^2-p)}{2(p-1)^2}=\frac{p(p+1)}{2}.$$

Let $n:=\frac{p(p+1)}{2}.$ Then from the assumption, we have $n\leq \frac{l-3}{2}.$
So by Theorem 5.17 (iii), we can choose  $n$ essentially distinct elements $c_l\Big(\frac{a_ix+b_i}{c_ix+d_i}\Big), 1\leq i \leq n,$ so that the cyclic subgroups $\Big\langle c_l\Big(\frac{a_ix+b_i}{c_ix+d_i}\Big) \Big\rangle$ are different, and that
$$ \bigcup_{t\in \{1,\, l-1\}\cup \mathfrak{Z}(l,p)}  GG_l(\mathbb{F}_p(x))^t\supseteq \bigcup_{i=1}^{n}\Big\langle c_l\Big(\frac{a_ix+b_i}{c_ix+d_i}\Big) \Big\rangle.$$

 Hence we get
$$\bigcup_{m=1}^{l-1} GG_l(\mathbb{F}_p(x))^{m}= \bigcup_{t\in \{1,\, l-1\}\cup \mathfrak{Z}(l,p)}  GG_l(\mathbb{F}_p(x))^t \supseteq\bigcup_{i=1}^{n}\Big\langle c_l\Big(\frac{a_ix+b_i}{c_ix+d_i}\Big) \Big\rangle.$$
This completes the proof.

 \hfill$\Box$

\bigskip

\bigskip

{\centerline{\bf 6. The Cases  {\bf ${5\leq l\leq 2n+1}$}
}}

\bigskip

\noindent Now, we consider the cases of $n> \frac{l-3}{2},$ i.e. $l\leq 2n+1,$ which seems difficult. For $n=2$ and $l=5,$ we have:

\bigskip

\noindent {\bf Theorem 6.1}\  {\it Assume that  $F$ is a field
 and $\Phi_5(x)$ is irreducible in $F[x]$.

  i) If ch$(F)=0,$ then $c(\mathfrak{T}_5(2; F))=4,$
 so $cs(\mathfrak{T}_5(2; F))=0.$

 ii) If ch$(F)=p\neq0, 2,$ then $c(\mathfrak{T}_{5}(2; F))=2(2+|\mathfrak{Z}(5,p)|).$ }

\noindent {\it Proof:} \ It suffices to prove
$$\beta= c_5(x)^{l_{1}}\cdot c_5(x+b)^{l_{2}}\notin G_5(F(x)),$$
where $b\neq0$ and $1\leq l_1, l_2\leq 4.$

Otherwise, if $\beta  \in G_5(F(x)),$ then in the proof of Theorem 5.17,
letting $n=2,$ we know that there exist two coprime polynomials
$f(x),g(x)\in F[x],$ with $f(x)$ monic, such that
$$\Phi_5(f,g) = \alpha  \Phi_5(x)^{e_1}\Phi_5(x+b)^{e_2}, \ \ \mbox{for some}\ \alpha\in F, \eqno(6.1)
$$
and that we have either deg$f=2$ or  deg$f=3.$
The proof of  case deg$f=3$
 is completely similar.

Now, we consider the case deg$ f=2.$
In this case, we have $e_1=e_2=1,$ so (6.1) becomes
$$\Phi_5(f,g) = \alpha  \Phi_5(x)\Phi_5(x+b).  \eqno(6.2)$$

Let $x=\zeta:=\zeta_5.$ Then we get $\frac{f(\zeta)}{g(\zeta)}=\zeta^i, 1\leq i \leq 4,$ so $f(\zeta)-\zeta^ig(\zeta)=0$ and $\zeta^{5-i}f(\zeta)-g(\zeta)=0.$
Hence we have
$$\Phi_5(x)\mid x^{2}f(x)-g(x)\ \ \mbox{or}\ \ \Phi_5(x)\mid f(x)-x^{2}g(x).$$
Similarly,   letting  $x=\zeta-b,$  we get
$$\Phi_5(x)\mid x^{2}f(x-b)-g(x-b)\ \ \mbox{or}\ \ \Phi_5(x)\mid f(x-b)-x^{2}g(x-b).$$

 Since $f(x)$ is monic, comparing the degrees   we have the following equalities
 $$\begin{array}{l}
                \Phi_5(x)=x^{2}f(x)-g(x) \ \ \mbox{or}\ \ -k_2\Phi_5(x)=f(x)-x^{2}g(x), \\
                \Phi_5(x)=x^{2}f(x-b)-g(x-b) \ \ \mbox{or}\ \ -k_2\Phi_5(x)=f(x-b)-x^{2}g(x-b),
            \end{array}
       \eqno(6.3)$$
where $k_2$ is the leading coefficient of $g(x).$

We claim that $k_2\neq0.$ Otherwise, if $k_2=0,$ then we have either
$$f(x)=x^{2}g(x)$$
or
$$f(x-b)=x^{2}g(x-b).$$

Let
$$f(x)=x^2+l_1x+l_0,\ \ \ g(x)=k_2x^2+k_1x+k_0.$$

If $f(x)=x^{2}g(x),$ then we have $f(x)=x^2, g(x)=1,$ so we get
$\Phi_5(x^2) = \Phi_5(x)\Phi_5(x+b). $
From $\Phi_5(x^2)=\Phi_5(x)\Phi_5(-x)$ we get
$\Phi_5(-x)=\Phi_5(x+b).$
Substituting $x=-\zeta$ we get $\Phi_5(b-\zeta)=0.$ Consequently $b-\zeta=\zeta^k$ for some $k=1,2,3,4.$ This is impossible,
since elements $1, \zeta$ and $1,\zeta,\zeta^k (k>1)$ are linearly independent over $F,$ since the minimal polynomial of $\zeta$ is of degree $4.$

If $f(x-b)=x^{2}g(x-b),$ then $f(x)=(x+b)^{2}g(x),$ so we get $f(x)=(x+b)^2, g(x)=1,$ therefore we have
$\Phi_5((x+b)^2) = \Phi_5(x)\Phi_5(x+b).$
Similarly, a contradiction arises.

Now, the formulas (6.3) lead to the following four cases:

(i). $\Phi_5(x)=x^{2}f(x)-g(x)=x^{2}f(x-b)-g(x-b).$

Hence $0\neq x^{2}(f(x)-f(x-b))=g(x)-g(x-b).$ This is impossible, since we have
$\mbox{deg}(g(x)-g(x-b))<\mbox{deg}g(x)\leq 2.$

 (ii). $-k_2\Phi_5(x)=f(x)-x^{2}g(x)=f(x-b)-x^{2}g(x-b).$

 Then $0\neq f(x)-f(x-b)=x^{2}(g(x)-g(x-b)).$
  This leads to a contradiction, since deg$f(x)=$deg$g(x)=2$ implies that
  $f(x)\neq f(x-b), g(x)\neq g(x-b), \mbox{deg}(f(x)-f(x-b))< 2.$

 (iii). $\Phi_5(x)=x^{2}f(x)-g(x)$ and $-k_2\Phi_5(x)=f(x-b)-x^{2}g(x-b).$

 From the first equality it follows that
 $f(x)=x^{2}+x+l_0, \ \ g(x)=(l_0-1)x^2-x-1.$
So the second equality gives
$-k_2=1-l_0=1-2b=2b(l_0-1)+1,$
so
$1-2b=-k_2=2b(l_0-1)+1$ $=2b(2b-1)+1.$
Since ch$(F)\neq2,$ we get $b=0,$ a contradiction.

 (iv). $\Phi_5(x)=x^{2}f(x-b)-g(x-b)$ and $-k_2\Phi_5(x)=f(x)-x^{2}g(x).$

From the second equality it follows that
$f(x)=x^{2}-k_2x-k_2,\ \  g(x)=k_2x^2+k_2x+k_2+1.$
 Then the first equality implies that
$2b+k_2=-1,\ \ 2bk_2-k_2=1.$
So $2b+2bk_2=0,$ therefore $k_2=-1.$ But this implies $b=0,$ a contradiction.

Thus, in all the four cases we get a contradiction.
In summary, the equality (6.2) does not hold.
\hfill$\Box$

\bigskip

\noindent {\it Remark} 6.2 \ \  The main result in [26] is  a special case of Theorem 6.1.
The assumption ch$(k)\neq 2$ is needed  in the proof of Theorem 6.1 because  for ch$(F)=2$ we have the equality:
$$\Phi_5(x^2+x, x^2+x+1)=\Phi_5(x)\Phi_5(x+1).$$

For $n=3$ and $l=5$ or $7,$ we have:
\bigskip

\noindent {\bf Theorem 6.3}\  {\it Assume that  $F$ is a field with ch$(F)\neq 2$ and that  $\Phi_l(x)$ is irreducible in $F[x],$ and assume that $l=5$ or $7.$

i) If ch$(F)=0,$ then $c(\mathfrak{T}_l(3; F))=6,$
so $cs(\mathfrak{T}_l(3; F))=0.$

ii) If  ch$(F)=p\neq 0,$ then $c(\mathfrak{T}_l(3;F))=3(2+|\mathfrak{Z}(l,p)|).$
}

\noindent {\it Proof:} \  Similar to the proof of Theorem 6.1, through a rather long computation, the proof can be achieved. \hfill$\Box$

\bigskip

\bigskip

{\centerline{\bf 7. Diophantine Equations
}}

\bigskip

\noindent To give a further example, we need the following two lemmas.

\bigskip

\noindent {\bf Lemma 7.1}\ {\it The  integer solutions of the diophantine equation
$$x^4+x^3y+x^2(y^2-1)+xy(y^2-1)+(y^2-1)^2=0$$
are only
$$(0,1), \ (-1,1), \ (0,-1), \ (1,-1).$$
In particular, if $y^2-1\neq0,$ then the equation has no integer solutions.}

\noindent {\it Proof:} \ Let $(x,y)=(a,b)$ be an integer solution.

If $b^2=1,$ then $b=\pm 1$ and $a^4+a^3b=0.$ It is easy to see that in these cases the solutions are only
$$(0,1), \ (-1,1), \ (0,-1), \ (1,-1).$$

If $b^2\neq 1,$ then rewrite the equation as
$$a^4=[(-ab)-(b^2-1)][a^2+b^2-1].$$
If $a=0,$ then $b^2-1=0,$ a contradiction; if $b=0,$ then $a^4-a^2+1=0,$ impossible. Hence, $ab\neq 0.$
Thus, we should have
$-ab>b^2-1> 0,$
so
$-ab\geq b^2.$
If $b> 0,$ then
$b\leq -a;$
if $b<0,$ then
$-b\leq a.$
So, in either cases, we have
$$a^4=[(-ab)-(b^2-1)][a^2+b^2-1]\leq[a^2-(b^2-1)][a^2+b^2-1]=a^4-(b^2-1)^2.$$
This is impossible since $b^2-1\neq0$. \hfill $\Box$

\bigskip

\noindent {\bf Lemma 7.2}\ {\it The  equation
$$x^4+x^3y+x^2(y^2+1)+xy(y^2+1)+(y^2+1)^2=0$$
has no real number solutions.}

\noindent {\it Proof:}
The polynomial can be written in the form
$$(x^4+x^3y+x^2y^2+xy^3+y^4)+(x^2+xy+y^2)+(y^2+1).$$
First two summands in brackets are nonnegative and the third is $\geq 1.$ Hence the value of the polynomial for $x,y \in \mathbb{R}$ is $\geq 1.$
\hfill$\Box$

\bigskip

\bigskip

{\centerline{\bf 8. A Further Example
}}

\bigskip

\noindent We continue to consider the cases of $l\leq 2n+1.$

We use the symbol $\mathfrak{S}^{*}_l(2; \mathbb{Z})$ denote a subgroup of $K_2(\mathbb{Q}(x))$
generated by   $2$ essentially distinct nontrivial elements of the form
$$c_l\left(\frac{a_1x+b_1}{c_1x+d_1}\right),\ \ c_l\left(\frac{a_2x+b_2}{c_2x+d_2}\right),$$
where
$$\left (\begin{matrix}
   a_1 & b_1\\
   c_1 & d_1
  \end{matrix}\right ),\ \ \left (\begin{matrix}
   a_2 & b_2\\
   c_2& d_2
  \end{matrix}\right )\in SL(2,\mathbb{Z})$$
  satisfying the `extra condition'
  $$\left (\begin{matrix}
   a_1 & b_1\\
   c_1 & d_1
  \end{matrix}\right )^{-1} \left (\begin{matrix}
   a_2 & b_2\\
   c_2& d_2
  \end{matrix}\right )\neq \pm\left (\begin{matrix}
   0 & -1\\
   1 & 1
  \end{matrix}\right ),\ \pm\left (\begin{matrix}
   1 & 0\\
   -1 & 1
  \end{matrix}\right ).$$

\bigskip

\noindent {\bf Theorem 8.1}\  {\it
  We have $c(\mathfrak{S}_5^{*}(2; \mathbb{Z}))=4,$
hence $cs(\mathfrak{S}_5^{*}(2; \mathbb{Z}))=0,$ i.e., $\mathfrak{S}_5^{*}(2; \mathbb{Z})$
contains no nontrivial cyclotomic subgroups.}

\noindent {\it Proof:}\ Let
$$\beta=c_5\left(\frac{a_1x+b_1}{c_1x+d_1}\right)^{l_{1}}\cdot c_5\left(\frac{a_2x+b_2}{c_2x+d_2}\right)^{l_{2}},$$
where
$$\left (\begin{matrix}
   a_1 & b_1\\
   c_1 & d_1
  \end{matrix}\right ), \left (\begin{matrix}
   a_2 & b_2\\
   c_2 & d_2
  \end{matrix}\right )\in SL(2,\mathbb{Z}).$$
We can assume $1\leq l_1,
l_2\leq 4.$

We claim that $\beta \notin G_5(\mathbb{Q}(x)).$ Otherwise, if $\beta  \in
G_5(\mathbb{Q}(x)),$ then as in the discussions of section 5, we know
that there exist two coprime polynomials $f(x),g(x)\in \mathbb{Q}[x]$   such
that
$$\Phi_5(f,g) = \alpha  \Phi_5(a_1x+b_1, c_1x+d_1)^{e_1}\Phi_5(a_2x+b_2, c_2x+d_2)^{e_2},  \ \mbox{where}\ \alpha \in \mathbb{Q}, \eqno(8.1)
$$
and that we have either deg$f=2$ or  deg$f=3.$

1. Case deg$ f=2.$

In this case, we have $e_1=e_2=1,$ so (8.1) becomes
$$\Phi_5(f,g) = \alpha \Phi_5(a_1x+b_1, c_1x+d_1)\Phi_5(a_2x+b_2, c_2x+d_2). $$

Let $X=\frac{a_1x+b_1}{c_1x+d_1}.$ Then, we have
$$\Phi_5\left((a_1-c_1X)^2f\left(\frac{d_1X-b_1}{a_1-c_1X}\right),(a_1-c_1X)^2g\left(\frac{d_1X-b_1}{a_1-c_1X}\right)\right)$$
 $$= \alpha \Phi_5(X)\Phi_5(a_2(d_1X-b_1)+b_2(a_1-c_1X), c_2(d_1X-b_1)+d_2(a_1-c_1X)). $$
So, it suffices to consider
$$\Phi_5(f,g) = \alpha  \Phi_5(x)\Phi_5(ax+b, cx+d),  \eqno(8.2)
$$
where $ad-bc=1$ and
   $$\left (\begin{matrix}
   a & b\\
   c& d
  \end{matrix}\right )\neq \pm\left (\begin{matrix}
   0 & -1\\
   1 & 1
  \end{matrix}\right ),\ \pm\left (\begin{matrix}
   1 & 0\\
   -1 & 1
  \end{matrix}\right ).$$

Noting that $\zeta \notin \mathbb{Q},$ by the action of the Galois group Gal$(\mathbb{Q}(\zeta)/\mathbb{Q}),$ we have
$$f(x)-\zeta g(x)=\alpha_{1}(x-\zeta^i)(ax+b-\zeta^j(cx+d)), \ \ \mbox{where} \ \ \alpha_1\in \mathbb{Q}(\zeta).
 \eqno(8.3)$$

Let
$$f(x)=x^2+l_1x+l_0,\ \ \ g(x)=k_2x^2+k_1x+k_0,$$
with $l_0,l_1,k_0,k_1,k_2\in \mathbb{Q}.$

Putting  this expressions in to (8.3) and comparing the coefficients,  we get
$$ck_2\zeta^{i+j+1}-c\zeta^{i+j}+(ck_1-dk_2)\zeta^{j+1}-ak_2\zeta^{i+1}+(d-cl_1)\zeta^j+a\zeta^i+(bk_2 -ak_1)\zeta+al_1-b=0,\eqno(8.4)$$
$$dk_2\zeta^{i+j+1}-d\zeta^{i+j}+ck_0\zeta^{j+1}-bk_2\zeta^{i+1}-cl_0\zeta^j+b\zeta^i-ak_0\zeta+al_0=0.\eqno(8.5)$$

We only consider the following  cases and the other cases are similar and easy.

1) If $i=1, j=2,$
from (8.4)(8.5) we have
$$ d-cl_1-ak_2=ck_2,\ al_1-b=ck_2, \eqno(8.6)$$
$$ck_1-dk_2-c=ck_2, \ a-ak_1+bk_2=ck_2, \eqno(8.7)$$
$$ck_0-d=dk_2,\ b-ak_0=dk_2, \eqno(8.8)$$
$$ cl_0+bk_2=-dk_2,\ al_0=dk_2. \eqno(8.9)$$

From $(8.6)$, we have
$(a^2+c^2+ac)k_2=1,$
so $k_2\neq 0;$
from $(8.7),$ we have
$ac+ c^2=-1;$
so from this equality and $(8.8)$,
we have
$(a+c)d=1-a^2;$
 therefore from
$(8.9)$, we have
$$1-a^2+ab=0.$$
Hence,  we have
$c^8+2c^6+4c^4+3c^2+1=0,$
 impossible.

2) If $i=1, j=3,$ then we have
$$ck_1-dk_2-c=-ak_2, \ a-ak_1+bk_2=-ak_2,$$
$$ d- cl_1=-ak_2, \ ck_2+al_1-b=-ak_2,$$
$$d-ck_0=bk_2,\ ak_0-b=bk_2,$$
$$ cl_0=bk_2,\ dk_2+al_0=-bk_2.,$$

From these equalities, we have $k_2\neq0$ and respectively
$$-1=(a^2+ac+c^2)k_2,\ \ \ a^2+ac-1=0,$$
 $$(ab+bc)k_2=1,\ \ \ c(b+d)=-ab.$$
Cancelling $k_2, b,c,$ we have
$a^8-2a^6+4a^4-3a^2+1=0,$
 impossible.

3) If $i=1, j=4,$ then we have
$$d-cl_1=ak_2=ck_2+a-ak_1+bk_2=al_1-b-c+ck_1-dk_2=0,$$
$$cl_0=bk_2=dk_2+b-ak_0=ck_0+al_0-d=0.$$
Then we have
$$a^2=a^2k_1,\ \ c^2k_1-cdk_2=c^2-1,\ \ c^2 k_0 - cd=0,\ \ b^2-abk_0=0.$$
Clearly $c\neq 0,$ and $a\neq0$ since $b^2-abk_0=0.$ So $k_1=1,$ therefore  from $c^2k_1-cdk_2=c^2-1,$ we have
$cdk_2=1;$ so from $ck_2+a-ak_1+bk_2=0,$ we have $b=-c.$
Hence from $b^2-abk_0=0, c^2k_0-cd=0,$ we get  $\frac{b}{a}=k_0=\frac{d}{c},$ that is, $ad-bc=0,$ a contradiction.

4) If $i=2, j=3,$ then we have
$$ck_1-dk_2=a,\ \ ck_2-ak_1+bk_2=a,\ \ d- ak_2-cl_1=a, \eqno(8.10)$$
$$al_1-b-c=a,\ \ ck_2=b,\ \ bk_2+cl_0=-b,\ \ dk_2-ak_0=b,\ \ al_0=b.\eqno(8.11)$$

From (8.10)(8.11), we have respectively
$$(c^2-1)k_2=a^2+ac, \ \ \ a^2k_2=1-a^2-ac-c^2.$$
So
$$a^4+ca^3+(c^2-1)a^2+(c^3-c)a+(c^2-1)^2=0.$$
From Lemma 7.1, we have $c^2=1.$ So $ a(a+c)=0.$

If $a=0,$ then $b=0$ from (8.11), a contradiction. So $a=-c,$ hence $a^2=1.$ From (8.10)(8.11), we have
$$ \left (\begin{matrix}
   a & b\\
   c& d
  \end{matrix}\right )= \left (\begin{matrix}
   a & 0\\
   -a & a
  \end{matrix}\right )= \pm\left (\begin{matrix}
   1 & 0\\
   -1 & 1
  \end{matrix}\right ).$$
This contradiction to the assumption.

5) If $i=j=2, $ we have
$$ck_1-dk_2-ak_2=-c,\ \ ak_1-bk_2=c,\eqno(8.12)$$
$$ a-cl_1+d=-c,\ \ al_1-b+ck_2=-c,\eqno(8.13)$$

From (8.16)(8.17), we have respectively
$$(a^2+1)k_2=c^2+ac,\ \ \ c^2k_2=-1-ac-a^2-c^2.$$
So  we have
$c^4+ac^3=-(a^2+1)(1+ac+a^2+c^2),$
i.e.,
$$c^4+ac^3+c^2(a^2+1)+ac(a^2+1)+(a^2+1)^2=0.$$
This contradicts  Lemma 7.2.

6) If $i=j=3,$ then we have
$$ck_1-dk_2-ak_2=ck_2,\ bk_2-ak_1-c=ck_2,\eqno(8.14)$$
$$a-cl_1+d=ck_2,\ al_1-b=ck_2,\eqno(8.15)$$
$$ck_0-bk_2=dk_2,\ ak_0+d=-dk_2,\eqno(8.16)$$

We claim that $k_2\neq0.$ In fact, if $k_2=0,$ then clearly $c=0$
(otherwise from (8.14)-(8.16), we will have $k_0=k_1=k_2=0$),  but from (8.15) this will
implies that $a=-d,$ impossible.

From (8.14)(8.15), we have respectively
$$(c^2+ac)k_2=a^2+1, \ \ (1+a^2)k_2=-(a^2+c^2+1). $$
So
$$c^4+ac^3+c^2(a^2+1)+ac(a^2+1)+(a^2+1)^2=0.$$
A contradiction arises from Lemma 7.2.

7) If $i=3, j=2,$ we have
$$ck_1-dk_2+a=-ak_2, \ \ \ d-cl_1=-ak_2,$$
$$ck_2-ak_1+bk_2=-ak_2,\ \ \ al_1-b-c=-ak_2,$$
$$b+ck_0=-bk_2,\ \ \ cl_0=bk_2,$$
$$ dk_2-ak_0=-bk_2,\ \ \ al_0-d=-bk_2.$$

From $cl_0=bk_2, al_0-d=-bk_2,$ we have
$cd=b(a+c)k_2;$
from $d-cl_1=-ak_2, $ $al_1-b-c=-ak_2,$ we have
$c^2-1=a(a+c)k_2.$
Hence $b=-c.$ Clearly $c\neq 0,$ otherwise, if $c=b=0,$ then $ad=1,$ hence $a^2=1,$ so from $ck_1-dk_2+a=-ak_2,$ we get $-k_2+a^2=-a^2k_2,$ that is, $1=0,$ a contradiction. Hence $c\neq 0.$

On the other hand, from $ck_1-dk_2+a=-ak_2, ck_2-ak_1+bk_2=-ak_2,$ we have
$$a^2=(1-a^2-ac-c^2)k_2.$$
In virtue of $c^2-1=a(a+c)k_2,$ we get
$$a^3(a+c)=(c^2-1)(1-a^2-ac-c^2).$$
So we obtain
$$a^4+a^3c+a^2(c^2-1)+ac(c^2-1)+(c^2-1)^2=0.$$
From Lemma 7.1, we get $c^2-1=0.$ Hence $ad=1-c^2=0.$

If $a\neq 0,$ then $d=0.$ So we have
$$k_0=1+k_2,\ l_0=-k_2,\ -ak_0=ck_2,\ al_0=ck_2.$$
Hence $k_0=-l_0=k_2=k_0-1,$ a contradiction.

Therefore $a=0.$ If $d=0,$ then
clearly $c(x)$ and $c(\frac{-c}{cx})$ are not essentially distinct, which contradicts the assumption.

Hence we get $a=0, d\neq 0.$ So we have
$$ck_1=dk_2, \ d=cl_1,\  k_0=1+k_2,\  -l_0=k_2,\  dk_2=ck_2,\  -d=ck_2.$$
Therefore $k_1=k_2\neq 0$ since $d\neq 0,$ so $d=c.$
Hence
 $$ \left (\begin{matrix}
   a & b\\
   c& d
  \end{matrix}\right )= \left (\begin{matrix}
   0 & -c\\
   c & c
  \end{matrix}\right )= \pm\left (\begin{matrix}
   0 & -1\\
   1 & 1
  \end{matrix}\right ).$$
This contradicts the assumption.

2. Case deg$f=3.$

In this case, we have $e_1+e_2=3,$ so by symmetry, it suffices to consider the case  $e_1=2, e_2=1,$ hence (8.1) becomes
$$\Phi_5(f,g) = \alpha \Phi_5(a_1x+b_1, c_1x+d_1)^2\Phi_5(a_2x+b_2, c_2x+d_2). $$

Similar as (8.2), it suffices to consider
$$\Phi_5(f,g) = \alpha  \Phi_5(x)^2\Phi_5(ax+b, cx+d),
$$
where $\left (\begin{matrix}
   a & b\\
   c & d
  \end{matrix}\right )\neq \left (\begin{matrix}
   0 & -1\\
   1 & 1
  \end{matrix}\right )\in SL(2,\mathbb{Z}).$
Similarly, we have
$$f(x)-\zeta g(x)=\alpha_{2}(x-\zeta^i)^2(ax+b-\zeta^j(cx+d)),  \ \ \mbox{where}\ \  \alpha_2\in \mathbb{Q}(\zeta).
\eqno(8.17)$$

Let
$$f(x)=x^3+l_2x^2+l_1x+l_0,\ \ \  g(x)=k_3x^3+k_2x^2+k_1x+k_0.$$
Putting these expression into (8.17) and comparing the coefficients, we have
$$2ck_3\zeta^{i+j+1}-2c\zeta^{i+j}-2ak_3\zeta^{i+1}+(ck_2-dk_3)\zeta^{j+1}+2a\zeta^i+(d-cl_2)\zeta^j$$
$$+(bk_3-ak_2)\zeta+(al_2-b)=0,$$
$$ck_3\zeta^{2i+j+1}-2dk_3\zeta^{i+j+1}-c\zeta^{2i+j}-ak_3\zeta^{2i+1}+2d\zeta^{i+j}+2bk_3\zeta^{i+1}-ck_1\zeta^{j+1}$$
$$+a\zeta^{2i}-2b\zeta^i+cl_1\zeta^{j}
+ak_1\zeta-al_1=0, $$
$$dk_3\zeta^{2i+j+1}-d\zeta^{2i+j}-bk_3\zeta^{2i+1}+b\zeta^{2i}-ck_0\zeta^{j+1}+cl_0\zeta^j+ak_0\zeta-al_0=0. $$

Similar to the proof of the case of deg$f=2,$ we can prove that these equalities do not hold. So we  omit the details of computations.

In summary, the equality (8.1) does not hold. So $\beta \notin G_5(\mathbb{Q}(x)),$ as claimed.
\hfill$\Box$

\bigskip

This example implies that the cases of $l\leq 2n+1$ are more complicated than imagination.

\bigskip

\noindent {\bf Question 8.2}\  {\it How to remove the condition $n\leq \frac{l-3}{2}$ in Theorem 5.17 ?}

\bigskip

\bigskip

{\centerline{\bf 9. The Cubes and Squares
}}

\bigskip

\noindent  From this section on, we will turn to the number field cases. In this section we will focus on the problem: When the cube or the square of a cyclotomic element is still cyclotomic ?  As a result,   we will construct some cyclotomic subgroups of order 5.

We need the following lemmas.

\bigskip

\noindent {\bf Lemma 9.1}(Selmer) \  {\it i) If $n\not \equiv 2  (\mbox{mod} \, 3),$ then the polynomials $x^n+x+1$ are irreducible in $\mathbb{Q}[x].$

ii) If $n \equiv 2  (\mbox{mod} \, 3),$ then the polynomials $x^n+x+1$ have a factor $x^2+x+1,$ but the polynomials $x^n+x+1/x^2+x+1$
 are still irreducible in $\mathbb{Q}[x].$ }

\noindent {\it Proof:} \ See [16]. \hfill$\Box$

\bigskip

\noindent {\bf Lemma 9.2} (Zsigmondy) \  {\it If $a>b>0,$ gcd$(a,b)=1$ and $n>1$ are positive integers, then $a^n+b^n$ has a prime factor that does not divide $a^k+b^k$ for all positive integers $k< n,$ with exception $2^3+1^3.$ }

\noindent {\it Proof:} \ See [30]. \hfill $\Box$

\bigskip

Let
$$f_{n,1}(x)= x^n+x+1,\ \ \mbox{if} \ \ n\equiv 1 \, (\mbox{mod} 3),$$
$$f_{n,2}(x)=x^n+x+1/x^2+x+1, \ \ \mbox{if}\ \ n\equiv 2 \, (\mbox{mod} 3).
$$

We can construct the cube of a cyclotomic element which is also cyclotomic as follows.

\bigskip

\noindent {\bf Theorem 9.3} \ {\it Assume that $p>3$ is a prime.  Let $\alpha$ be a zero of $f_{p,i}(x),$ where $i=1$ or $2,$ and $F=\mathbb{Q}(\alpha).$  Then we have
$$1\neq c_p(\alpha)^{3}= c_p(\alpha^3) \in G_p(F).$$}
\noindent {\it Proof:} \ Clearly $\alpha^{p-1}\neq 0, 1.$ From $\alpha^p+\alpha+1=0,$ we have
$$\alpha^{p-1}(\alpha^p+\alpha+1)=0,$$
therefore
$$\alpha^{2p-1}+\alpha^{p-1}=\alpha+1,$$
so
$$\alpha^{2p}+\alpha^p+1=\alpha^2+\alpha+1,$$
that is,
$$\frac{\alpha^{3p}-1}{\alpha^p-1}=\frac{\alpha^3-1}{\alpha-1},$$
which implies
$$\frac{\alpha^p-1}{\alpha-1}=\frac{\alpha^{3p}-1}{\alpha^3-1}.$$
Hence
$$\Phi_p(\alpha)=\Phi_p(\alpha^3),$$
and therefore
 $$c_p(\alpha)^3=\{\alpha^3, \Phi_p(\alpha)\}=\{\alpha^3, \Phi_p(\alpha^3)\}=c_p(\alpha^3) \in G_p(F).$$

Now, we  prove that $c_p(\alpha)^3\neq 1.$ Since $p> 3$ is a prime, it suffices to prove $c_p(\alpha)\neq 1.$

 At first, we can simplify the formula for $c_p(\alpha).$ Namely,
 $$\Phi_p(\alpha)=\frac{1-\alpha^p}{1-\alpha}=\frac{1+(\alpha+1)}{1-\alpha}=\frac{\alpha+2}{1-\alpha}.$$
 Hence
 $$c_p(\alpha)=\{\alpha, \Phi_p(\alpha)\}=\left\{\alpha, \frac{\alpha+2}{1-\alpha}\right\}=\{\alpha,\alpha+2\},$$
 since $\{\alpha, 1-\alpha\}=1.$ But
 $$\{\alpha,\alpha+2\}=\left\{-2\left(\frac{-\alpha}{2}\right), 2\left(1+\frac{\alpha}{2}\right)\right\}
 =\left\{-2, 1+\frac{\alpha}{2}\right\}\left\{\frac{-\alpha}{2}, 2\right\}=\{-2,2+\alpha\}\{\alpha, 2\}.$$
So we have
$$c_p(\alpha)=\{-2,2+\alpha\}\{\alpha, 2\}.$$

 Clearly, $\alpha$ is a unit. Hence $v_{\mathfrak{p}}(\alpha)=0$ for every prime ideal $\mathfrak{p}.$ Therefore, for every prime ideal $\mathfrak{p}\nmid 2,$ we get
$$\tau_{\mathfrak{p}}(c_p(\alpha))= \tau_{\mathfrak{p}}(\{-2,\alpha+ 2\}\{\alpha, 2\})\equiv (-2)^{v_{\mathfrak{p}}(\alpha+ 2)}(\mbox{mod}\, \mathfrak{p}).  \eqno(9.1)$$

When $p\equiv 1 (\mbox{mod}\, 3),$  from Lemma 9.1,  $f_{p,1}(x)$ is irreducible in $\mathbb{Q}[x].$
So the minimal polynomial of $\alpha+2$ is
$$f_{p,1}(x-2)=(x-2)^p+(x-2)+1=x^p-2px^{p-1}+\ldots +2^{p-1}px+x-(2^p+1).$$
Hence $N_{F/\mathbb{Q}}(\alpha+2)=2^p+1.$

When $p\equiv 2 (\mbox{mod}\, 3),$  from Lemma 9.1,  $f_{p,2}(x)$ is also irreducible in $\mathbb{Q}[x].$
So the minimal polynomial of $\alpha+2$ is $f_{p,2}(x-2).$ From
$$(x-2)^p+(x-2)+1=[(x-2)^2+(x-2)+1]f_{p,2}(x-2).$$
We know that $N_{F/\mathbb{Q}}(\alpha+2)=\frac{1}{3}(2^p+1).$

Now we can assume that $p\equiv 1 (\mbox{mod}\, 3)$ since the case of $p\equiv 2 (\mbox{mod}\, 3)$ can be treated in similar way.

Suppose that we have the decomposition of prime ideals
 $$(\alpha+2)\mathcal{O}_F=\mathfrak{p}_1^{e_1}\mathfrak{p}_2^{e_2}\cdots \mathfrak{p}_m^{e_m}.$$
 In virtue of  $N_{F/\mathbb{Q}}(\alpha+2)=2^p+1$, we can assume that $e_i\geq 1$ and $m\geq 1.$ Let $p_i$ be primes (need not  be different) such that
 $(p_i)=\mathfrak{p}_i\cap \mathbb{Z}.$ Then
 $$N_{F/\mathbb{Q}}((\alpha+2)\mathcal{O}_F)=N_{F/\mathbb{Q}}(\mathfrak{p}_1)^{e_1}\cdots N_{F/\mathbb{Q}}(\mathfrak{p}_m)^{e_m}={p_1}^{e_1f_1}\cdots {p_m}^{e_mf_m}\mathbb{Z},$$
 where $f_i=f(\mathfrak{p}_i | p_i)$ are the residue class degrees.

 From Lemma 9.2, the number $2^p+1$ has a primitive prime divisor, say, $q,$ i.e., $q \mid2^p+1$ but $q\nmid 2^d+1$ for any integer $1\leq d < p.$

 Assume that  $v_q(2^p+1)=l.$ Then we have
 $$N_{F/\mathbb{Q}}((\alpha+2)\mathcal{O}_F)=N_{F/\mathbb{Q}}(\alpha+2)\mathbb{Z}=(2^p+1)\mathbb{Z}=q^la\mathbb{Z},$$
 where $q\nmid a.$ Therefore
 $$q^la\mathbb{Z}={p_1}^{e_1f_1}\cdots {p_m}^{e_mf_m}\mathbb{Z}.$$
 This implies that  $q$ must be one  of the primes $p_1, p_2, \ldots, p_m,$ say $q=p_1.$ Note that the primes $p_i$ may not be distinct. So we have
 $$l=e_1f_1+\ldots \geq e_1.$$

On the other hand, clearly we have $q\neq 3,$ i.e., $q\geq 5,$ so
$$5^l\leq q^l< 2^p+1,$$
therefore $l< p$, hence $e_1\leq l< p,$ that is,  $v_{\mathfrak{p}_1}(\alpha+2)=e_1< p.$ This implies that
$$p_1=q \nmid 2^{v_{\mathfrak{p}_1}(\alpha+2)}+1,$$

 On the other hand, we also have
 $$p_1\nmid  2^{v_{\mathfrak{p}_1}(\alpha+2)}-1.$$
 In fact, otherwise if $p_1|  2^{v_{\mathfrak{p}_1}(\alpha+2)}-1,$ then from $p_1|2^p+1$ we get
$$p_1| (2^p+1)+ (2^{v_{\mathfrak{p}_1}(\alpha+2)}-1)=2^p+ 2^{v_{\mathfrak{p}_1}(\alpha+2)}=2^{v_{\mathfrak{p}_1}(\alpha+2)}(2^{p-v_{\mathfrak{p}_1}(\alpha+2)}+1).$$
 Since $p_1=q \neq 2,$ we have $p_1|2^{p-v_{\mathfrak{p}_1}(\alpha+2)}+1.$ This contradicts the choice of $q=p_1$ since $v_{\mathfrak{p}_1}(\alpha+2)\neq 0.$

Hence from (9.1), we get
$$\tau_{\mathfrak{p}_1}(c_p(\alpha))\equiv (-2)^{v_{\mathfrak{p}_1}(\alpha+ 2)} \not\equiv 1 (\mbox{mod}\, \mathfrak{p}_1).$$
Therefore $c_p(\alpha)\neq 1.$
\hfill$\Box$

\bigskip

\noindent {\bf Lemma 9.4} \ {\it The Galois group of the polynomial $x^p+x+1$ is isomorphic to  $S_p,$ the symmetric group of degree $p$.}

\noindent {\it Proof:} It follows from Lemma 9.1 and Theorem 1 in  [11].
\hfill$\Box$

\bigskip

\noindent {\bf Lemma 9.5} \ {\it Let $L/F$ be a Galois extension of finite degree $n$ with Galois group $G:=$ Gal$(L/F).$ Then the kernel of the canonical homomorphism $K_2(F)\rightarrow K_2(L)^G$ is killed by $n$.}

\noindent {\it Proof:} See [19].
\hfill$\Box$

\bigskip

\noindent {\bf Corollary 9.6} \ {\it The assumption is the same as in Theorem 9.3. Let $\widetilde{F}$ be the normal closure of $F=\mathbb{Q}(\alpha)$ with Galois group Gal$(\widetilde{F}/\mathbb{Q}).$  Then}

(i) {\it For any $\sigma \in\mbox{Gal}(\widetilde{F}/\mathbb{Q}),$ we have}
$$c_p(\sigma(\alpha))^{\pm 3}\neq 1.$$

(ii) {\it In $K_2(\widetilde{F})$ we have the equality:
$$\prod_{\sigma \in G} c_p(\sigma(\alpha))=c_p(-2),$$
i.e. the element $\prod_{\sigma \in G} c_p(\sigma(\alpha))$ is also cyclotomic.}

\noindent {\it Proof:} \ (i) From Lemma 9.4, we have $[\widetilde{F}:F]=[\widetilde{F}:\mathbb{Q}]/[F:\mathbb{Q}]=|S_p|/p=(p-1) !$ and  from Lemma 9.5, we know that  the kernel of the homomorphism  $K_2(F)\rightarrow K_2(\widetilde{F})^G \subseteq K_2(\widetilde{F})$ is killed by $[\widetilde{F}:F].$ By $([\widetilde{F}:F], p)=1,$ we get the injection:
 $$G_p(F)\hookrightarrow  G_p(\widetilde{F}),$$
 since $G_p(F)$ is contained in the $p$-torsion of $K_2(F)$ (see [1]).

 Then the result follows from Theorem 9.3 and the facts $\sigma (c_p( \alpha))=c_p(\sigma(\alpha))$ and $c_p( \alpha)^{-1}=c_p( \alpha^{-1}).$

(ii) From the proof of Theorem 9.3, in $K_2(\widetilde{F})$ we have:
$$\prod_{\sigma \in G} c_p(\sigma(\alpha))=\prod_{\sigma \in G} \{-2, 2+\sigma(\alpha)\}\{\sigma(\alpha),2\}\ \ \ \ \ \ \ \ \ \ \ \ \ \ $$
$$\ \ \ \ \ \ \ \ \ \ \ \ \ \ =\{-2, \prod_{\sigma \in G}(2+\sigma(\alpha))  \}\cdot \{\prod_{\sigma \in G}\sigma(\alpha), 2\}$$
$$\ \ \ \ \ \ \ \ \ \ \ \ \ \ \ \ \ \ \ \ \ \ \ \ \ \ \ =\{-2, -f_{p,1}(-2)\}\{N_{\widetilde{F}/\mathbb{Q}}(\alpha), 2\}\ \ \ \ \ \ \ \ \ \ \ \ \ \ \ \ \ \ \ \ $$
$$\ \ \ \ \ \ \ \ =\{-2, 2^{p}+1\} =\Big\{-2, \frac{2^{p}+1}{3}\Big\}$$
$$=\{-2, \Phi_p(-2)\}=c_p(-2).$$
\hfill$\Box$

In  the case of the polynomial $f_{p,2}(x),$ Browkin told me that $\langle c_p(\alpha)\rangle$ is a cyclotomic subgroup when $p=5.$ Moreover, in the following example, we can construct three different nontrivial cyclotomic subgroups.

\bigskip

\noindent {\bf Example 9.7}\   \ Let $p=5.$ Then it is easy to show that
$$f_{5, 2}(x)=x^3-x^2+1.$$

Let $\alpha$ be a zero of $f_{5,2}(x)$ and $F=\mathbb{Q}(\alpha).$ Then  from Theorem 9.3,  we get
  $$1\neq c_5(\alpha)^{3}= c_5(\alpha^3) \in G_5(F).$$

On the other hand, we have
$$c_5(\alpha)^2=c_5(\alpha)^{-3}=c_5(\alpha^{-3}),$$
$$c_5(\alpha)^{4}=c_5(\alpha)^{-1}=c_5(\alpha^{-1}).$$
Hence we conclude that $\langle c_5(\alpha)\rangle\subset G_5(F),$ i.e., $\langle c_5(\alpha)\rangle$ is a cyclotomic subgroup.

Moreover, let $F_{n}=F(\sqrt[5^{n-1}]{\alpha}\,).$ By the formula $\Phi_{5^n}(x)=\Phi_{5}(x^{5^{n-1}}),$ we get
$$c_{5^n}(\sqrt[5^{n-1}]{\alpha}\,)^{5^{n-1}}=\{\sqrt[5^{n-1}]{\alpha}, \Phi_{5^n}(\sqrt[5^{n-1}]{\alpha}\,)\}^{5^{n-1}}=\{\alpha, \Phi_{5}(\alpha)\}=c_5(\alpha).$$
So $G_{5^n}(F_n)$ also contains the cyclotomic subgroup $\langle c_{5^n}(\sqrt[5^{n-1}]{\alpha} \,)^{5^{n-1}}\rangle= \langle c_5(\alpha)\rangle.$

Now, let $\widetilde{F}$ be the normal closure of $F=\mathbb{Q}(\alpha).$ Then from Theorem 2 in [11], we have Gal$(\widetilde{F}/\mathbb{Q})\cong S_3,$ so $[\widetilde{F}:F]=|S_3|/[F:\mathbb{Q}]=2.$ Thus from Theorem 9.5, we have the injection:
$$G_5(F)\hookrightarrow G_5(\widetilde{F}).$$
Hence for any $\sigma \in G:=$Gal$(\widetilde{F}/\mathbb{Q}),$ we have
$$1\neq c_5(\sigma(\alpha))^{\pm 3}\in G_5(\widetilde{F}).$$
and therefore
$$\bigcup_{\sigma\in G}\langle c_5(\sigma(\alpha)) \rangle  \subseteq G_5(\widetilde{F}). \eqno(9.2)$$

Let $\alpha:=\alpha_1, \alpha_2$ and $\alpha_3$ be the three roots of the polynomial $f_{5, 2}(x)=x^3-x^2+1.$
Then (9.2) becomes
$$\langle c_5( \alpha_1) \rangle\cup \langle c_5(\alpha_2) \rangle \cup \langle c_5(\alpha_3) \rangle \subseteq G_5(\widetilde{F}).$$

\noindent {\bf Claim:}\ {\it The  cyclotomic subgroups
$\langle c_5( \alpha_1) \rangle,  \langle c_5(\alpha_2) \rangle ,  \langle c_5(\alpha_3) \rangle$ are different
 from each other. Hence $G_5(\widetilde{F})$ contains   at least three nontrivial cyclotomic subgroups.}

In fact, from the proof of Theorem 9.3, we have
$$c_5(\alpha_1)=\{-2, 2+\alpha_1\}\{\alpha_1,2\},$$
and $N_{F/\mathbb{Q}}(2+\alpha_1)=11.$ So as in the proof of Theorem 9.3, we can prove that $(2+\alpha_1)$ is a prime ideal in $\mathcal{O}_F.$

Since $[\widetilde{F}:F]=2,$ we have
$$(2+\alpha_1)\mathcal{O}_{\widetilde{F}}=\mathfrak{p}_1, \mathfrak{p}_1^2 \ \mbox{or}\ \mathfrak{p}_1\mathfrak{p}_2,$$
 where $  \mathfrak{p}_1, \mathfrak{p}_2$ are prime ideals of $\mathcal{O}_{\widetilde{F}}.$ Note that $\alpha$ is a unit. Then we have
$$\tau_{\mathfrak{p}_1}(c_5(\alpha_1))\equiv  -2 \, (\mbox{mod} \, \mathfrak{p}_1),\ \ \mbox{if}\ (2+\alpha_1)\mathcal{O}_{\widetilde{F}}=\mathfrak{p}_1,\ \ \ \ \ \ \ \ \ \ \ \ \ \ \ $$
$$\tau_{\mathfrak{p}_1}(c_5(\alpha_1))\equiv  4 \, (\mbox{mod} \, \mathfrak{p}_1),\ \ \mbox{if}\ (2+\alpha_1)\mathcal{O}_{\widetilde{F}}=\mathfrak{p}_1^2,\ \ \ \ \ \ \ \ \ \ \ \ \ \ \ \ \ \eqno(9.3)$$
$$\ \ \ \ \ \ \ \ \ \ \ \ \ \ \ \ \ \ \ \ \ \ \ \ \ \tau_{\mathfrak{p}}(c_5(\alpha_1))\equiv  -2 \, (\mbox{mod} \, \mathfrak{p}),\ \ \mbox{if}\ (2+\alpha_1)\mathcal{O}_{\widetilde{F}}=\mathfrak{p}_1\mathfrak{p}_2, \ \mbox{and}\ \mathfrak{p}=\mathfrak{p}_1\ \mbox{or}\ \mathfrak{p}_2.\ \ \ \ \ \ \ \ \ \ \ \ \ \ \ \ $$
We can do similar works for the elements $\alpha_2$ and $\alpha_3.$

Assume that  $\langle c_5( \alpha_1) \rangle=\langle c_5(\alpha_2) \rangle.$ Then  we have $  c_5( \alpha_1) =  c_5(\alpha_2)^i, $ for some $1\leq i \leq 4.$
 But the prime ideals over $2+\alpha_1, 2+\alpha_2$ and $2+\alpha_3$ are different from each other. So
 we have
 $$\tau_{\mathfrak{p}}(c_5(\alpha_2))\equiv 1  \, (\mbox{mod} \, \mathfrak{p}),$$
 where $\mathfrak{p}$ is the prime ideal appearing in (9.3).
 Hence from (9.3) we have:
 $$(-2)^i\ \mbox{or} \ 4^i \equiv 1  \, (\mbox{mod} \, \mathfrak{p}),\  \mbox{for some} \ 1\leq i \leq 4.$$
  But $\mathfrak{p} $ is over $11,$ so we get
 $$(-2)^i\ \mbox{or} \ 4^i \equiv 1  \, (\mbox{mod} \, 11),\ \mbox{for some}\ 1\leq i \leq 4.$$
It is easy to check that  this is impossible. Hence we must have $\langle c_5( \alpha_1) \rangle \neq  \langle c_5(\alpha_2) \rangle. $ Similarly we can prove that
$\langle c_5( \alpha_1) \rangle \neq  \langle c_5(\alpha_3) \rangle $ and $\langle c_5( \alpha_2) \rangle \neq  \langle c_5(\alpha_3) \rangle. $
The claim is proved.

\bigskip

 \noindent {\it Question:}\   How many cyclotomic subgroups are there in $G_5(\widetilde{F}) $ ?

\bigskip

 Furthermore, from the proof of Theorem 9.3 and similar to the proof of Corollary 9.6 (ii), we can show   the following equality   in $K_2(\widetilde{F}):$
$$c_5(\alpha_1)c_5(\alpha_2)c_5(\alpha_3)= \{-2,11\}=\{-2, \Phi_5(-2)\}.  $$

By Lemma 9.5, we have the injection $K_2(\mathbb{Q})\hookrightarrow K_2(\widetilde{F}).$ While in $K_2(\mathbb{Q})$ the tame symbol of $\{-2,11\}$ is
$$\tau_{11}(\{-2,11\})\equiv -2 \not \equiv 1\, (\mbox{mod}\, 11),$$
so we get
$$c_5(\alpha_1)c_5(\alpha_2)c_5(\alpha_3)=c_5(-2)\neq 1.$$
  \hfill$\Box$

\bigskip

We can also construct  a cyclotomic element in $K_2$ of some quadratic field such that its square is also a cyclotomic element.

\bigskip

\noindent {\bf Lemma 9.8} \ {\it For any integer $n\geq 1$ and any prime $p,$ the polynomials $f(x)=x^n+x^{n-1}+p$ are  irreducible over $\mathbb{Q}.$}

\noindent {\it Proof:} \ Clearly we can assume that $n\geq 2.$ The Newton polygon of $f(x)$ for the prime $p$ has vertices $(0,1), (n-1,0),(n,0).$ Therefore this polygon has two sides with slopes $1/(n-1)$ and $0,$ respectively.

It follows that in $\mathbb{Q}_p[x]$ this is reducible $f(x)=f_1(x)f_2(x),$ where deg$f_1=n-1,$ deg$f_2=1.$ Any root of $f_1(x)$ generates an extension of $\mathbb{Q}_p$ of degree $n-1,$ by the value of the corresponding slope. Consequently $f_1(x)$ is irreducible in $\mathbb{Q}_p[x].$

Consequently, if $f(x)$ were reducible in $\mathbb{Q}[x],$ then it should have factors of degrees $1$ and $n-1.$ It is impossible since $f(x)$ does not vanish at $\pm 1, \pm 2,$ so it does not have a root in $\mathbb{Q}.$

Thus $f(x)$ is irreducible in $\mathbb{Q}[x].$ \hfill$\Box$

\bigskip

\noindent {\it Remark} 9.9  \ i) We can also give a more computational proof of Lemma 9.8 as follows (see [11]).

Assume that we have the decomposition
$$x^n+x^{n-1}+ p=f(x)g(x),\ \ \mbox{where\  deg}f(x), \mbox{deg}g(x)\geq 1.$$
 Since $p$ is a prime, we can assume that the constant term of, say, $f(x)$ is $\pm 1.$

If $f(x)$ has a root of unity $\alpha,$ that is, $\alpha^n+\alpha^{n-1}+ p=0$ with $|\alpha|=1,$
then
$$p=|\alpha^n+\alpha^{n-1}|=|\alpha^{n-1}||\alpha+1|=|\alpha+1|.$$
Clearly $|\alpha+1|< 2$ if $\alpha\neq 1.$ So $\alpha=1.$ But $1$ is not a root of $x^n+x^{n-1}+p,$ a contradiction.

Hence,  $f(x)$ has no roots of unity. This implies deg$f(x)\geq 2$ and  $f(x)$ must have a root $\alpha$ with $|\alpha|< 1$.
So we have
$$p=|\alpha^n+\alpha^{n-1}|\leq |\alpha^n|+|\alpha^{n-1}|<2,$$
a contradiction again.
These contradictions prove the irreducibility of the polynomial $x^n+x^{n-1}+ p.$

ii) Similarly, we can prove that $x^n+x^{n-1}-p$ is also irreducible if $n\geq 1$ and $p\geq 3$ is a prime.

\bigskip

\noindent {\bf Theorem 9.10}\ {\it  Assume that  $p\geq 3$ be a prime. Let $\alpha$ be a zero of the polynomial $x^{p}+x^{p-1}+2$ and $F=\mathbb{Q}(\alpha).$
Then we have
$$1\neq c_p(\alpha)^{2}= c_p(\alpha^2) \in G_p(F).$$}
\noindent {\it Proof:} \ From $\alpha^{p}+\alpha^{p-1} + 2=0,$ we have
$$(1+ \alpha )^2\Phi_p(-\alpha)=(1+\alpha)(\alpha^p+1)=\alpha(\alpha^p+\alpha^{p-1})+1+\alpha=1-\alpha.$$
Then from $\{\alpha, (1+\alpha)^2\}=\{-1, 1+\alpha \}^2\{-\alpha, 1+\alpha \}^2=\{1, 1+\alpha \}=1,$ we get
$$c_p(\alpha)=\{\alpha, \Phi_p(\alpha)\}=\{\alpha, (1-\alpha)\Phi_5(\alpha)\}=\{\alpha, (1+\alpha)^2\Phi_p(-\alpha)\Phi_p(\alpha)\}=\{\alpha, \Phi_5(\alpha^2)\}.$$
So
$$c_p(\alpha)^2=\{\alpha^2, \Phi_p(\alpha^2)\}=c_p(\alpha^2)\in G_p(F).$$

Now, we prove $c_p(\alpha)\neq 1.$

From $\alpha^p+\alpha^{p-1}+2=0$ and $(1+ \alpha )^2\Phi_p(-\alpha)=1-\alpha,$ we have
$$\Phi_p(\alpha)=\frac{\alpha^p-1}{\alpha-1}=\frac{\alpha^p+1}{\alpha+1}\cdot \frac{\alpha+1}{\alpha-1}+\frac{2}{1-\alpha}=\frac{1-\alpha}{(1+\alpha)^2}\cdot \frac{1+\alpha}{\alpha-1}+\frac{2}{1-\alpha}=\frac{1+3\alpha}{1-\alpha^2}.$$
So
$$c_p(\alpha)=\{\alpha, \Phi_p(\alpha)\}=\left\{\alpha, \frac{1+3\alpha}{1-\alpha^2}\right\}=\{-3, 1+3\alpha\}\{-1, 1+\alpha\}^{-1}.$$

From Lemma 9.8, $x^p+x^{p-1}+2=0$ is irreducible over $\mathbb{Q}.$ So we have
$$N_{F/\mathbb{Q}}(1+3\alpha)=2(3^p+1).$$

Similar to the proof of Theorem 9.3, we can choose a primitive prime factor of $3^p+1,$ say $q.$
Clearly $q\neq 2,$ otherwise we would have $q|3+1,$ which  contradicts the choice of $q.$

Let $\mathfrak{p}$ be a prime lying above $q.$ Then similarly we can show that
$$1\leq v_{\mathfrak{p}}(1+3\alpha)\leq v_q(3^p+1)< p.$$

From $\mathfrak{p}|1+3\alpha=(1+\alpha)+2\alpha$ and $\mathfrak{p}|q\neq 2,$ we know that $\mathfrak{p}\nmid 1+\alpha,$ i.e., $v_{\mathfrak{p}}(1+\alpha)=0.$
So
$$\tau_{\mathfrak{p}}(c_p(\alpha))=
                (-3)^{v_{\mathfrak{p}}(1+3\alpha)}(-1)^{v_{\mathfrak{p}}(1+\alpha)}\equiv (-3)^{v_{\mathfrak{p}}(1+3\alpha)} \not \equiv 1(\mbox{mod}\, \mathfrak{p}).$$
Therefore $c_p(\alpha)\neq 1.$ \hfill$\Box$

\bigskip

\noindent {\bf Example 9.11}\ When $p=5,$ let $\alpha$ be a zero of the polynomial $x^{5}+x^{4}+2$ and $F=\mathbb{Q}(\alpha).$ Similar to Example 9.7, we conclude that $\langle c_5(\alpha)\rangle$ is a cyclotomic subgroup  of order $5$. Note that when $p=3,$ $\langle c_3(\alpha)\rangle$ is a cyclotomic subgroup of order $3.$

\bigskip

We can also construct a quadratic field $F$ such that $G_5(F)$ contains a cyclotomic subgroup of order $5$.
This was suggested to me by Browkin.

\bigskip

\noindent {\bf Example 9.12}\  The roots of the polynomial $x^2-3x+1=0$ are $\frac{3\pm \sqrt{5}}{2}.$ Let $\beta=\frac{3+\sqrt{5}}{2},$ and $F=\mathbb{Q}(\beta)=\mathbb{Q}(\sqrt{5}).$ Then we have
$$\Phi_5(-\beta)=(1-\beta^2)^2.$$
In virtue of  $\{\beta,(1-\beta^2)^2\}=\{\beta^2,1-\beta^2\}=1,$ we get
$$c_5(\beta)=\{\beta, (1-\beta^2)^2\Phi_5(\beta)\}=\{\beta, \Phi_5(-\beta)\Phi_5(\beta)\}=\{\beta, \Phi_5(\beta^2)\}.$$
So
$$c_5(\beta)^2=\{\beta^2, \Phi_5(\beta^2)\}=c_5(\beta^2)\in G_5(F).$$

  Similar as Example 9.7, we know that $\langle c_5(\beta)\rangle$ is a nontrivial cyclotomic subgroup.
But we need to prove that $c_5(\beta)\neq 1.$

In fact, note that $\beta^2+1=3\beta.$ Then we have
$$\Phi_5(\beta)=(1+\beta^2)^2-\beta^2+\beta(1+\beta^2)=9\beta^2-\beta^2+3\beta^2=11\beta^2.$$
Consequently,
$c_5(\beta)=\{\beta,11\beta^2\}=\{\beta, 11\}.$

In $\mathcal{O}_F=\mathbb{Z}[(1+\sqrt{5})/2],$ we have
$11=(4+\sqrt{5})(4-\sqrt{5}).$ Therefore $4+\sqrt{5}$ generates a prime $\mathfrak{p}$.
From $\beta^2-3\beta +1=0,$ we get $\beta(3-\beta)=1$ and $(\beta-1)^2=\beta,$ so $(3-\beta)(\beta-1)^2=1.$ These imply that $\beta, \beta-1$ are both units.
So we have $v_{\mathfrak{p}}(11)=1$ and $v_{\mathfrak{p}}(\beta)=0,$ and therefore
$$\tau_{\mathfrak{p}}(c_5(\beta))=\tau_{\mathfrak{p}}(\{\beta, 11\})\equiv \beta \not\equiv 1 (\mbox{mod} \mathfrak{p}).$$

Moreover, let  $\overline{\beta}=\frac{3-\sqrt{5}}{2}.$ Then similarly we have $c_5(\overline{\beta})= \{\overline{\beta}, 11\},$
and it is easy to see that  $c_5(\beta)c_5( \overline{\beta})=1.$ So we get
$\langle c_5(  \beta)\rangle  = \langle c_5( \overline{\beta})\rangle.$

\bigskip

\noindent {\it Question:} \ Are there any  nontrivial cyclotomic subgroups other than $\langle c_5(  \beta)\rangle$ which is contained in $G_5(\mathbb{Q}(\sqrt{5}))$ ?

\bigskip

 We do not know how to construct other cyclotomic subgroups. In particular, we do not know whether $\langle c_7(\alpha)\rangle$ is a cyclotomic subgroup with $\alpha$ as described above.

\bigskip

\bigskip

{\centerline{\bf 10. The Non-Closeness
}}

\bigskip

\noindent  In this section, for any number field $F,$ we will construct a subgroup generated by an infinite number of cyclotomic elements to the power of some prime,  which contain no nontrivial cyclotomic elements. This is more clear than what Browkin's conjecture implies.

We need the following celebrated  result.

\bigskip

\noindent {\bf Theorem 10.1}\ (Faltings [3])\  {\it Any smooth, projective curve over a number field $F$ that has genus greater than $1$ can have only finitely many $F$-rational points.} \hfill$\Box$

\bigskip

In the following, we will use the symbols $g(C)$ and $g(F(C))$ to denote respectively the genus of a curve $C$ and its function field $F(C).$
We also need a genus formula on Kummer extensions of function fields.

 Let $K/k$ be an algebraic function field where $k$ is the field of constants and contains a primitive $m$-th root of unity (with $m>1$ and $m$ relatively prime to the characteristic of $k$). Suppose that $u\in K$ is an element satisfying
$$u\neq w^d\ \ \mbox{for all}\ \ w\in K\ \ \mbox{and}\ \ d|m, d>1.$$
Let
$$K^{\prime}=K(y)\ \ \mbox{with}\ \ y^m=u.$$
Such an extension $K^{\prime}/K$ is said to be a Kummer extension of $K.$ We have the following genus formula.

\bigskip

\noindent {\bf Lemma 10.2}\ ([17]) \ {\it Let $K^{\prime}/K$ be the Kummer extension of function field $K$ with $y^m=u$ as above.
If $k^{\prime}$ denotes the constant field of $K^{\prime},$ then
$$g(K^{\prime})=1+\frac{m}{[k^{\prime}:k]}\left(g(K)-1+\frac{1}{2}\sum_{P\in S_K}\left(1-\frac{r_P}{m}\right )deg P \right)$$
where $r_P:=gcd(m, v_P(u))$ and $S_K$ is the set of places of $K/k.$} \hfill $\Box$

\bigskip

\noindent {\bf Lemma 10.3} \ {\it Let $F$ be a number field. Assume that $n\geq 3$ and $p$ is a prime. If either $p\geq 5$
 or $p=2$ but $n\neq 3,4,5, 6, 8, 10, 12$ or $p=3$ but $n\neq 3,4, 6,$  then there are only finitely many $F$-rational points on the curve $C: \Phi_n(x)=cy^p,$ where $c\in F^*.$}

{\it Proof.} \ Let $\overline{C}$ be the projective closure of $C$ over $F$ i.e.
 $$\overline{C}: \Phi_n(x, z)-c y^p z^{\varphi(n)-p}=0.$$

Note that $\overline{C}$ is a singular curve with singular point $(0:1:0).$  So we need to consider the normalization of $\overline{C},$ i.e.,
$$\pi: {\overline{C}}^{\prime}\longrightarrow \overline{C}.$$
As we know([5]), ${\overline{C}}^{\prime}$ is a projective smooth curve over $F.$ It is also well known that the genus of a projective smooth curve is equal to the genus of its function field([6]). So we have
$$g({\overline{C}}^{\prime})=g(F({\overline{C}}^{\prime})).$$
Since $\pi$ is a birational morphism, we have $F({\overline{C}}^{\prime})\simeq F(\overline{C})\simeq F(C),$ so
$g(F(\overline{C}^\prime))=g(F(C)),$ therefore
$$g({\overline{C}}^{\prime})=g(F(C)).$$

Now, we calculate the genus $g(F(C)).$

At first, since $F$ is a perfect field, the genus is unchanged under the algebraic extension of $F.$ So $g(F(C))=g(\overline{F}(C)),$ where $\overline{F}$ is the algebraic closure of $F.$

Clearly, we have
$$\overline{F}(C)=\overline{F}(x, y)=\overline{F}(x)(y)\ \ \mbox{ with} \ \ y^p=\Phi_n(x).$$
It is easy to see that $\overline{F}(x)(y)/\overline{F}(x)$ is a Kummer extension. As is well-known, the genus of the rational function field $\overline{F}(x)$ is trivial, i.e., $g(\overline{F}(x))=0.$

For the Kummer extension  $\overline{F}(x)(y)/\overline{F}(x)$ with
$$ y^p=u:=\Phi_n(x)=\prod_{\begin{array}{ll}
1\leq i \leq n,
(n,i)=1
 \end{array}}(x-\zeta^i),$$
where $\zeta$ is the $n$-th primitive root of unity,  it is easy to show that for any $P\in S_{\overline{F}(x)},$ we have

i) if $P=(x-\zeta^i), 1\leq i \leq n,$ gcd$(n,i)=1,$ then $v_{P}(u)=1,$ so $r_{P}=$gcd$(p, v_{P}(u))=1;$

ii) if $P=(x-a), a\neq \zeta^i,$ gcd$(n,i)=1, $ then $v_{P}(u)=0,$ so $r_{P}=p;$

iii) if $P=\infty=(\frac{1}{x}),$ then $v_{\infty}(u)=-\varphi(n),$ so $r_{\infty}=\mbox{gcd}(p, \varphi(n)).$

We apply Lemma 10.2 to the extension $\overline{F}(x)(y)/\overline{F}(x).$ Note that $\overline{F}$ is an algebraically closed field, so the constant field of $\overline{F}(x)(y)$ is also $\overline{F}$ and so deg$P=1$ for any place $P\in S_{\overline{F}(x)}.$ Therefore we get
$$g(\overline{F}(C))=1+p\left[-1+\frac{1}{2}\varphi(n)\left(1-\frac{1}{p}\right)+\frac{1}{2}\left(1-\frac{\mbox{gcd}(p,\varphi(n))}{p}\right)\right].$$
Thus, to prove $g(\overline{F}(C))\geq 2,$ it suffices to prove
$$\varphi(n)(p-1)> p+\mbox{gcd}(p, \varphi(n)). \eqno(10.1)$$

Note that  we have $\varphi(n)\geq 2$ since $n\geq 3.$

For $p\geq 5,$ if $\varphi(n)\geq 3,$ then $\varphi(n)(p-1)\geq 3(p-1)> 2p\geq p+\mbox{gcd}(p, \varphi(n));$ if $\varphi(n)=2,$ then
$\varphi(n)(p-1)=2(p-1)>p+1= p+\mbox{gcd}(p, \varphi(n)).$

For $p=2,$ the inequality (10.1) becomes
$$\varphi(n)> 2+\mbox{gcd}(2, \varphi(n)). $$
It is easy to see that this inequality holds if and only if $\varphi(n)>4.$ So $n\neq 3, 4, 5, 6,8,10,12.$

For $p=3,$ the inequality (10.1) becomes
$$2\varphi(n)> 3+\mbox{gcd}(3, \varphi(n)).  $$
Obviously, this holds if and only if $\varphi(n)>3.$ So $n\neq 3,4,6.$

Summarily,  we have  $g(F(C))\geq 2$ under the assumption on $n$ and $p.$
 So $g(\overline{C}^{\prime})\geq 2.$
Hence, $\overline{C}^{\prime}$ is a projective smooth curve of genus $\geq 2.$ Therefore, from  Theorem 10.1,  there are only finitely many $F$-rational points on $\overline{C}^{\prime},$ while $\pi$ is an $F$-birational morphism, hence there are also only finitely many $F$-rational point on $\overline{C}$ and therefore on $C,$ as required.

 \hfill$\Box$

\bigskip

\noindent {\bf Theorem 10.4} \ {\it Assume that $F$ is a number field and
 $n\neq 1,4,8,12$ is a positive integer. If there is a prime $p$ such that  $p^2|n,$
then  there exist infinitely many  nontrivial cyclotomic elements $\alpha_1, \alpha_2, \ldots, \alpha_m, \ldots  \in
G_{n}(F)$ so that
 $$\langle\alpha_1^{p}\rangle\subsetneq \langle\alpha_1^{p}, \alpha_2^p\rangle\subsetneq \ldots \subsetneq \langle\alpha_1^{p}, \alpha_2^p, \ldots, \alpha_m^p\rangle \subsetneq \ldots $$
 and
$$\langle\alpha_1^{p}, \alpha_2^p, \ldots, \alpha_m^p, \ldots \ \rangle\cap G_{n}(F)=\{1\}.$$
}

\noindent {\it Proof:} \ Let $S$ be a finite set of places of $F$ containing all archimedean ones, and all places above $p$ and above the primes ramified in $F.$ Moreover, we assume that $S$ is sufficiently large, so that the ring $\mathcal{O}_{F,S}$ of $S$-integers is a unique factorization domain. Let $P_S$ denote the set of all the rational primes which the finite primes in $S$ lie above.

Let $\mathbb{J}=\{1,2,\ldots,\frac{n}{p}-1\},$  and let
$N$ be a positive integer which is greater than $p,$ the rational
primes ramified in $F$ and all the rational primes in $P_S$.

Note that the polynomials $\Phi_n(x)$ and $\Phi_n^{\prime}(x)$ are coprime, so there exist two polynomials $g(x), h(x)\in \mathbb{Z}[x]$ and an integer $m_0$ so that
$$g(x)\Phi_n(x)+h(x)\Phi_n^{\prime}(x)=m_0. \eqno(10.2)$$

Let $M_1=m_0\prod_{1<q\leq N}q$ with $q$ running over all the rational primes less than $N$.
We can choose a sufficiently large integer $k_1$ and a rational
prime $p_1$ such that $p_1\mid \Phi_{n}(k_1M_1)$ (so $p_1 \nmid k_1M_1$).

Let
$$A_1:=\left\{\begin{array}{ll}
                k_1M_1, &  \text{if} \ v_{p_1}(\Phi_{n}(k_1M_1))=1,\\
                k_1M_1+p_1, & \text{if} \ v_{p_1}(\Phi_{n}(k_1M_1))> 1.
            \end{array}
        \right.$$
Then it is easy to show that $v_{p_1}(\Phi_{n}(A_1))=1,$ i.e., $p_1\parallel \Phi_{n}(A_1).$ In fact, if $v_{p_1}(\Phi_{n}(k_1M_1))> 1,$ then from Taylor formula, we
$$\Phi_{n}(k_1M_1+p_1)=\Phi_{n}(k_1M_1)+\Phi^{\prime}_{n}(k_1M_1)p_1+\frac{1}{2}\Phi^{\prime\prime}_{n}(k_1M_1)p_1^2+\ldots .$$
We must have $p_1 \nmid \Phi^{\prime}_{n}(k_1M_1).$ Otherwise, if $p_1 | \Phi^{\prime}_{n}(k_1M_1),$ then from (10.2), we have $p_1 | m_0.$ But according to the choice of $M_1,$ we have $m_0 | M_1,$ so $p_1| k_1M_1,$ a contradiction. Therefore we have $v_{p_1}(\Phi_{n}(A_1))=v_{p_1}(\Phi_{n}(k_1M_1+p_1)=1,$ as claimed.

Let
$$M_2=\prod_{q_1\mid k_1M_1}q_1\prod_{q^{\prime}_1\mid k_1M_1+p_1}q^\prime_1\prod_{q_2\mid\Phi_{n}(k_1M_1)}q_2\prod_{q^{\prime}_2\mid\Phi_{n}(k_1M_1+p_1)}q^\prime_2,$$
where $q_1,q_2$ run over rational primes. Then we can choose a
sufficiently large integer $k_2$ and a rational prime $p_2$ such
that $p_2\mid \Phi_{n}(k_2M_2),$ and similarly we get $A_2$ with $p_2 \parallel \Phi_n(A_2).$ 

Repeating this procedure, we get
the following sequences of elements of $K_2(F)$:
$$\{c_n(A_i)^{pj}\mid i=1,2,\ldots\},\ \ j\in \mathbb{J}, \eqno(10.3)$$
where
$$A_i=\left\{\begin{array}{ll}
                k_iM_i, &  \text{if} \ v_{p_i}(\Phi_{n}(k_iM_i))=1,\\
                k_iM_i+p_i, & \text{if} \ v_{p_i}(\Phi_{n}(k_iM_i))> 1,
            \end{array}
        \right.$$
in which $p_i$ is a rational prime satisfying $p_i\mid
\Phi_{n}(k_iM_i)$ (therefore $p_i\nmid k_iM_{i}$) and
$$M_i=\prod_{q_1\mid k_iM_{i-1}}q_1\prod_{q_1\mid  k_iM_{i-1}+p_{i-1}}q^\prime_1\prod_{q_2\mid\Phi_{n}(k_iM_{i-1})}q_2\prod_{q_2\mid\Phi_{n}(k_iM_{i-1}+p_{i-1})}q^\prime_2.$$
Hence $p_i\parallel \Phi_{n}(A_i).$ Note that we have $p_i\notin P_S$ for any $i.$

\bigskip

\noindent {\bf Claim 1} \ For each $j\in \mathbb{J},$ the elements of (10.3) are all nontrivial and different
from each other.

\bigskip

In fact, for each $p_i,$ we can choose a prime $\mathfrak{p}_i \subset \mathcal{O}_{F,S}$ with
$\mathfrak{p}_i| p_i$ since $p_i\notin P_S.$ According to the above construction, $p_i$ is
unramified in $F,$ so from $p_i\| \Phi_{n}(A_i)$ and
$\mathfrak{p}_i\mid p_i,$ we have $\mathfrak{p}_i
\| \Phi_{n}(A_i),$ i.e.,
$v_{\mathfrak{p}_i}(\Phi_{n}(A_i))=1.$ So
$$\tau_{\mathfrak{p}_i}(c_n(A_i)^{pj})\equiv A_i^{pj}(\mbox{mod}\, \mathfrak{p}_i).$$

It suffices to prove $A_i^{pj}\not \equiv 1(\mbox{mod}\, \mathfrak{p}_i).$ In fact, otherwise, 
assume that $A_i^{pj}\equiv 1(\mbox{mod}\, \mathfrak{p}_i).$ Let

$$j=p^{m}j_1, \ \ \mbox{where} \ \ 0\leq m\leq v_p(n)-2\  \mbox{and} \ (j_1, p)=1.$$
Then
$$\mbox{gcd}(n, p j )=p^{m+1}\cdot \mbox{gcd}\Big(\frac{n}{p^{m+1}}, j_1\Big)=p^{m+1}\cdot \mbox{gcd}\Big(\frac{n}{p^{v_p(n)}}, j_1\Big)$$
since $p\nmid j_1.$

So from $A_i^{n}\equiv 1(\mbox{mod}\, \mathfrak{p}_i)$
 and $A_i^{pj}\equiv 1(\mbox{mod}\, \mathfrak{p}_i),$ we have
 $$A_i^{p^{m+1}\cdot gcd( np^{-v_p(n)},j_1)}\equiv 1(\mbox{mod}\, \mathfrak{p}_i).$$
  Therefore $A_i^{\frac{n}{p}}\equiv 1(\mbox{mod}\, \mathfrak{p}_i).$

 It is easy to prove that
 there exists a polynomial $\Psi_{n,p}(x)\in \mathbb{Z}[x]$ such that
 $$\Phi_{n}(x)\Psi_{n,p}(x)=\Phi_p(x^{\frac{n}{p}}).$$

  Hence we get
$$0\equiv\Phi_{n}(A_i)\Psi_{n,p}(A_{i})=\Phi_p(A_i^{\frac{n}{p}})\equiv p \, (\mbox{mod}\, \mathfrak{p}_i),$$
that is, $\mathfrak{p}_i \mid p.$ This is impossible since $p_i\neq p$ (note that $p\in P_S$). So we get
$$\tau_{\mathfrak{p}_i}(c_n(A_i)^{pj})\equiv A_i^{pj}\not\equiv 1(\mbox{mod}\, \mathfrak{p}_i),$$
 which implies that
$c_n(A_i)^{pj}$ is nontrivial.

Next, we have $p_{i+1}\nmid M_{i+1},$ so $\mathfrak{p}_{i+1}\nmid M_{i+1}.$ Therefore, according to the construction, $\mathfrak{p}_{i+1}\nmid A_l, \mathfrak{p}_{i+1}\nmid \Phi_{n}(A_{l}), l\leq i, $ hence
$$\tau_{\mathfrak{p}_{i+1}}(c_n(A_l)^{pj})\equiv 1 (\mbox{mod}\, \mathfrak{p}_{i+1}),\ \ \forall \ l\leq i. $$
But from the above discussion, we know that
$$\tau_{\mathfrak{p}_{i+1}}(c_n(A_{i+1})^{pj})\not \equiv 1 (\mbox{mod}\, \mathfrak{p}_{i+1}). $$
Hence
$$c_n(A_l)^{pj}\neq c_n(A_{i+1})^{pj}, \ \ \forall \, l\leq i.$$
The claim is proved.

\bigskip

\noindent {\bf Claim 2}\ There exist some $i_1$  so that $c_n(A_{i_1})^{pj}\notin G_{n}(F)$ for each $j\in \mathbb{J}.$

\bigskip

At first, if there are only finitely many $i$ such that $c_n(A_i)^p\in G_n(F),$ choose a large integer $N_1$ so that
when $i>N_1,$
 $c_n(A_i)^p\notin G_n(F);$ otherwise,  we can choose an infinite subset $I_1\subseteq \mathbb{N}$
 so that for any $i\in I_1,$ we have  $c_n(A_i)^p\in G_n(F).$

  Next, if
there are only finitely many $i\in I_1$ such that $c_n(A_i)^{2p}\in G_n(F),$ choose a large integer $N_2>N_1$
(if $N_1$ exists) so that when $i\in I_1$ and $i>N_2,$ we have
 $c_n(A_i)^{2p}\notin G_n(F);$ otherwise, we choose an infinite subset $I_2\subseteq I_1$
 so that for any $i\in I_2,$ we have  $c_n(A_i)^{2p}\in G_n(F).$

Repeating this procedure, finally  we will get  an infinite set $I\subseteq \mathbb{N}$
and a set of integers:
 $$J:=\{j_1,j_2,\ldots, j_s\}, \ \ \mbox{with} \ \ 1\leq
j_1<j_2<\ldots<j_s\leq \frac{n}{p}-1,$$
 which satisfy
$$c_n(A_{i})^{pj}\in G_{n}(F), \  \ i\in I,\ \ j\in J$$
and
$$c_n(A_{i})^{pj}\notin G_{n}(F),  \  \ i\in I, \ \  j\in \mathbb{J}- J.$$

In the above construction, if $J=\emptyset,$ i.e., if for each $j\in \mathbb{J},$ it is always the first case, in another words, there are  only finitely many $i$ such that $c_n(A_i)^{pj}\in G_n(F),$ then the proof of the claim  has been done. Otherwise, we have $J\neq \phi,$ that is, $J$ is nonempty. We will prove that this is impossible.

In fact,  since $c_n(A_{i})^{pj}\in G_{n}(F),$ for $i\in I,\ \ j\in J,$ we can assume that
$$c_n(A_{i})^{pj}=c_n(B_{ij}), \ \ \mbox{where}\ \ i\in I,\ j\in J, B_{ij}\in F^{*}.$$

By the Dilichlet-Hasse-Chevalley theorem (see [21]), the group of $S$-units in $\mathcal{O}_{F,S}$ is finitely generated: There are
fundamental $S$-units $\varepsilon_{1}, \varepsilon_{2}, \ldots,\varepsilon_t$ such that every $S$-unit can be written in the form
$$\zeta^{r}\varepsilon_{1}^{k_1}\varepsilon_{2}^{k_2}\cdots \varepsilon_{t}^{k_{t}}, \ \ \mbox{where}\ r,k_1,\ldots,k_{t}\in \mathbb{Z}.$$
Here $\zeta$ is a generator of the group of roots of unity in $F$ and $0\leq r <$ ord$\zeta.$

By Lemma 10.3, the
equation $\Phi_{n}(x)=cy^p$ has only finitely many solutions with
$x,y\in F.$ Hence, there are only finitely many $x\in F$ such that
$\Phi_{n}(x)$ can be written in the form of $cy^p$ with $c$ having the form:
$$\zeta^{r}\varepsilon_{1}^{k_1}\varepsilon_{2}^{k_2}\cdots \varepsilon_{t}^{k_{t}},\ \ \ \ 0\leq r< p, 0\leq k_j< p,1\leq j<t. \eqno(10.4)$$

 Hence, we can find an integer $\widetilde{N}\in I$ such that when $i\in I$ and $i> \widetilde{N},$ $\Phi_{n}(B_{ij})$ can not be written in the form of $cy^p,$ where $c$ has the form of (10.4).  This implies that  we must have
 $$\Phi_{n}(B_{ij})=c_{ij}a_{ij}y^{p}_{ij},$$
 where $c_{ij}$ has the form of (10.4), and $a_{ij}\in F^*\backslash \mathcal{O}^*_{F,S}\cdot(F^*)^p$ and $ y_{ij}\in F^*.$

Assume that
$$a_{ij}\mathcal{O}_{F,S}=\mathfrak{q}_{ij1}^{e_{ij1}}\mathfrak{q}_{ij2}^{e_{ij2}}\cdots \mathfrak{q}_{ijs}^{e_{ijs}}.$$

 We claim that there must exists some $k_0$, where $1\leq k_0\leq s,$ such that $p\nmid e_{ijk_0},$ i.e., $p\nmid v_{\mathfrak{q}_{ijk_0}}(a_{ij}).$

In fact, if $p\mid e_{ijk}$ for $1\leq k\leq s,$ letting $e_{ijk}=pe_{ijk}^\prime, 1\leq k\leq s,$   we have
$$a_{ij}\mathcal{O}_{F,S}=(\mathfrak{q}_{ij1}^{e^\prime_{ij1}}\mathfrak{q}_{ij2}^{e^\prime_{ij2}}\cdots \mathfrak{q}_{ijs}^{e^\prime_{ijs}})^p=(a^\prime_{ij}\mathcal{O}_{F,S})^p,\ \ \mbox{for some}\ a^\prime_{ij}\in F^*,$$
since, according to the choice of $S,$ $\mathcal{O}_{F,S}$ is a UFD, so a PID.  Therefore
$$a_{ij}=u_{ij}(a^{\prime \prime}_{ij})^p,\ \ \mbox{where}\ u_{ij}\in \mathcal{O}^*_{F,S}  \ \mbox{and}\ a^{\prime \prime}_{ij}\in F^*,$$
that is, $a_{ij}\in \mathcal{O}^*_{F,S}\cdot(F^*)^p,$ a contradiction. So the claim is true.

 For the convenience, we denote $\mathfrak{q}_{ij}:=\mathfrak{q}_{ijk_0}.$

Therefore, from Claim 3, we conclude that if $i> \widetilde{N},$ then for each $j\in J$ there must exists a
prime  $\mathfrak{q}_{ij}$ such that
$$p\nmid v_{\mathfrak{q}_{ij}}(\Phi_{n}(B_{ij})).$$
Since $p\in P_S,$ we have $\mathfrak{q}_{ij}\nmid p.$

Now, we prove that this will lead to a contradiction. 

On one hand,
we have
$$c_n(B_{ij})^{\frac{n}{p}}=c_n(A_{i})^{nj}=1.$$
On the other hand,
if $v_{\mathfrak{q}_{ij}}(B_{ij})>0,$ then $v_{\mathfrak{q}_{ij}}(\Phi_{n}(B_{ij}))=0,$ a contradiction; if $v_{\mathfrak{q}_{ij}}(B_{ij})<0,$ then from $v_{\mathfrak{q}_{ij}}(\Phi_{n}(B_{ij}))=v_{\mathfrak{q}_{ij}}(B_{ij})\cdot$deg$\Phi_{n}(x)$ and $p\mid$ deg$\Phi_{n}(x),$ we have
$$p\mid v_{\mathfrak{q}_{ij}}(\Phi_{n}(B_{ij})),$$
a contradiction again. Hence, we must have $v_{\mathfrak{q}_{ij}}(B_{ij})=0.$

Note that $v_{\mathfrak{q}_{ij}}(\Phi_{n}(B_{ij}))>0,$ i.e., $\mathfrak{q}_{ij}\mid \Phi_{n}(B_{ij}).$
Computing the tame symbol, we get
$$\tau_{\mathfrak{q}_{ij}}(c_n(B_{ij})^{\frac{n}{p}})\equiv B_{ij}^{v_{\mathfrak{q}_{ij}}(\Phi_{n}(B_{ij}))\frac{n}{p}}(\mbox{mod}\, \mathfrak{q}_{ij}).$$
In virtue of $c_n(B_{ij})^{\frac{n}{p}}=1,$ we have
$$B_{ij}^{v_{\mathfrak{q}_{ij}}(\Phi_{n}(B_{ij}))\frac{n}{p}}\equiv 1(\mbox{mod}\, \mathfrak{q}_{ij}).$$
From $\mathfrak{q}_{ij}\mid \Phi_{n}(B_{ij}) \mid (B_{ij}^{n}-1),$ we obtain  $B_{ij}^{n}\equiv 1(\mbox{mod}\, \mathfrak{q}_{ij})$.
Hence, we get
$$B_{ij}^{\frac{n}{p}}\equiv 1(\mbox{mod}\, \mathfrak{q}_{ij}),$$
 since  gcd$(n, v_{\mathfrak{q}_{ij}}(\Phi_{n}(B_{ij}))\frac{n}{p})=\frac{n}{p}.$
Therefore, we conclude
$$0\equiv \Phi_{n}(B_{ij})\Psi_{n,p}(B_{ij})=\Phi_{p}(B_{ij}^{\frac{n}{p}})\equiv \Phi_{p}(1)\equiv p \,(\mbox{mod}\, \mathfrak{q}_{ij}),$$
i.e., $\mathfrak{q}_{ij}\mid p,$ a contradiction. Thus,  Claim 2 is proved.

Now,
let $\alpha_1=c_n(A_{i_1}).$ Then Claim 2 implies
$$\langle \alpha_1^{p}\rangle \cap G_n(F)=\{1\},$$
as required.

Next, we construct $\alpha_2.$

At first, from Claim 1, we can choose a sufficient large integer $N_1$ so that when $i>N_1+i_1,$  we have $c_n(A_i)^{pj}\not\in \langle \alpha_1\rangle$ for any $j\in \mathbb{J}.$

Let
$$M^\prime_i:=M_{N_1+i_1+i}, A^\prime_i:=A_{N_1+i_1+i},p_i^\prime:=p_{N_1+i_1+i}.$$
Th notations are the same   as above.

As (7.4),  we  construct  sequences of elements:
$$\{c_n(A_{i_1})^{pj}\cdot c_n(A^{\prime}_i)^{pj^{\prime}}\mid i=1,2,\ldots\},\ \ j, j^{\prime} \in \mathbb{J}, \eqno(10.5)$$
with $p^{\prime}_i\parallel \Phi_{n}(A^{\prime}_i).$

Similarly, as Claim 1, we can prove that for fixed each pair $(j,j^\prime) \in \mathbb{J}\times \mathbb{J},$ the elements of (10.5) are all nontrivial and different from each other.

Assume that for each $i,$ there exists a couple $(j,j^\prime)\in
\mathbb{J}\times \mathbb{J}$ such that
$$c_n(A_{i_1})^{pj}\cdot c_n(A^{\prime}_i)^{pj^{\prime}}\in
G_n(F).$$

Similar to the above discussion, there exists an infinite subset
$I^\prime\subseteq \mathbb{N}$ and $J^\prime \subseteq \mathbb{J}\times \mathbb{J}$ such
that
$$c_n(A_{i_1})^{pj}\cdot c_n(A^{\prime}_i)^{pj^{\prime}}\in G_n(F), \ \ i\in I^\prime,\ (j,j^\prime)\in J^\prime$$
and
$$c_n(A_{i_1})^{pj}\cdot c_n(A^{\prime}_i)^{pj^{\prime}}\notin G_n(F), \ \ i\in I^\prime,\ (j,j^\prime)\in \mathbb{J}\times \mathbb{J}-J^\prime.$$

Now, assume that
$$c_n(A_{i_1})^{pj}\cdot c_n(A^{\prime}_i)^{pj^{\prime}}=c_n(B^\prime_{ij}), \ \ i\in I^\prime,\ (j,j^\prime)\in J^\prime,$$
with $B^\prime_{ij}\in F^{*}.$ As above, we can prove similarly that
$J'=\emptyset.$

Hence, there must exist some $i_2$  so that
$$c_n(A_{i_1})^{pj}\cdot c_n(A_{i_2})^{pj^{\prime}}\notin
G_{n}(F), \ \ \mbox{for any}\ \ (j,j^\prime)\in \mathbb{J}\times \mathbb{J}.$$

 Let
$\alpha_2=c_n(A_{i_2}).$ Since $i_2>N_1+i_1,$ we have $\alpha_2\notin \langle \alpha_1\rangle$. So we get
$$\langle\alpha_1^{p}\rangle\subsetneq \langle\alpha_1^{p}, \alpha_2^p\rangle$$
$$\langle \alpha_1^{p}, \alpha_2^p\rangle \cap G_n(F)=\{1\}.$$

Repeating the procedure,   we can find  $\alpha_1,
\alpha_2, \ldots, \alpha_m, \ldots  \in G_{n}(F)$ so that
 $$\langle\alpha_1^{p}\rangle\subsetneq \langle\alpha_1^{p}, \alpha_2^p\rangle\subsetneq \ldots \subsetneq \langle\alpha_1^{p}, \alpha_2^p, \ldots, \alpha_m^p\rangle\subsetneq \ldots $$
 and
$$\langle\alpha_1^{p}, \alpha_2^p, \ldots, \alpha_m^p, \ldots \ \rangle\cap G_{n}(F)=\{1\}.$$
The proof is finished.$\Box$

\bigskip

\noindent {\bf Conjecture 10.5}\ {\it Let $F$ be a number field. If $p>5$ is a prime, then $G_p(F)$ contains no nontrivial cyclotomic subgroups.}  $\Box$

\bigskip

\noindent\textbf{Acknowledgement} We are grateful to Professor Jerzy Browkin for his many helpful suggestions which make the proofs of the results in this paper more transparent.

\bigskip

\noindent {\large K}{\footnotesize EJIAN} {\large X}{\footnotesize U}\ \  \ \ \ \ \ $\verb"kejianxu@amss.ac.cn"$

\noindent College of Mathematics

\noindent Qingdao University

\noindent  Qingdao 266071

\noindent China

\bigskip

\noindent {\large C}{\footnotesize HAOCHAO} {\large S}{\footnotesize UN}\ \  \ \ \ \ \ $\verb"sunuso@163.com"$

\noindent School of Mathematics

\noindent Jilin University

\noindent Changchun 130012

\noindent China

\end{document}